\documentclass[12pt,a4paper]{article}

\usepackage[utf8]{inputenc}
\usepackage{amsmath,amssymb,amsthm}
\usepackage{mathrsfs}
\usepackage{hyperref}
\usepackage{geometry}
\usepackage{xcolor}
\usepackage{indentfirst}

\newcommand{\corr}[1]{\textcolor{blue}{#1}}
\numberwithin{equation}{section}
\usepackage{enumitem}

\newtheorem{theorem}{Theorem}[section]
\newtheorem{lemma}[theorem]{Lemma}
\newtheorem{remark}[theorem]{Remark}
\newtheorem{proposition}[theorem]{Proposition}

\newenvironment{proofof}[1]
{\par\noindent\textbf{Proof of #1.}}
{\qed\par}
\hypersetup{
	pdftitle={Ground States of a Coupled Fractional Schrodinger System under a Total-Mass Constraint},
	pdfauthor={Lijun Liu, Jiabin Zuo},
	colorlinks = true,
	linkcolor  = black,   
	urlcolor   = black,   
	citecolor  = black,   
	breaklinks = true    
}

\newcommand{\R}{\mathbb{R}}

\newcommand{\Hd}{H^s(\R^N)}
\newcommand{\fl}{(-\Delta)^s}

\newcommand{\norm}[2]{\left\|#1\right\|_{#2}}
\newcommand{\abs}[1]{\left|#1\right|}
\newcommand{\intR}{\int_{\R^N}}

\title{Ground States of a Coupled Fractional Schr\"{o}dinger System under a Total-Mass Constraint: Complete Classification and Concentration Behavior}
\author{
	Lijun Liu$^{a}$, Jiabin Zuo$^{a,*}$\\[6pt]
	{\small $^{a}$\ School of Mathematics and Information Science, Guangzhou University,}\\[1pt]
	{\small Guangzhou 510006, China}\\[2pt]
	{\small $^{*}$\ Corresponding author.}\\[2pt]
	{\small E-mail addresses: \href{mailto:jun9059@e.gzhu.edu.cn}{jun9059@e.gzhu.edu.cn} (L.~Liu), \href{mailto:zuojiabin88@163.com}{zuojiabin88@gzhu.edu.cn} (J.~Zuo).}
}
\date{}

\begin{document}

\maketitle

\begin{abstract}
	We consider a coupled fractional Schr\"{o}dinger system with trapping potentials under a total-mass constraint, in which both the intraspecies and the interspecies interactions are attractive. For every fractional order \(s\in(0,1)\) and dimension \(N\ge1\), we obtain a complete classification of the existence and non-existence of minimizers: minimizers exist if and only if both intraspecies coefficients lie below the single-component critical mass and the interspecies coupling coefficient lies below an explicit critical threshold $\beta^*$, and no minimizer exists otherwise. The threshold $\beta^*$ is characterized by means of an auxiliary intensity function, which yields a lower bound for the energy functional without invoking a new coupled Gagliardo--Nirenberg inequality. As the coupling coefficient approaches $\beta^*$, we analyze the concentration behavior. In the one-dimensional half-Laplacian case, the two components are shown to blow up at a common zero of the trapping potentials and, after rescaling, to converge to the same fractional ground-state profile, with the mass ratio determined explicitly by the intraspecies coefficients and the critical mass. For polynomial trapping potentials whose flattest common zero has order $q_0\in(0,1)$, the sharp blow-up rate is established, and both components concentrate simultaneously at the flattest common minimum.
\end{abstract}

\vspace{6pt}

\vspace{4pt}
\section{Introduction}

	In this paper, we are interested in the following coupled fractional Schr\"{o}dinger system with trapping potentials in $\mathbb{R}^N$:
\begin{equation}\label{eq:sys}
	\begin{cases}
		(-\Delta)^s u_1 + V_1(x)u_1 = \mu u_1 + a_1|u_1|^{p-2}u_1 + \beta |u_1|^{q-2}|u_2|^q u_1, \\
		(-\Delta)^s u_2 + V_2(x)u_2 = \mu u_2 + a_2|u_2|^{p-2}u_2 + \beta |u_1|^q |u_2|^{q-2} u_2,
	\end{cases}
\end{equation}
where $0<s<1$, $N \ge 1$, $p = 2 + \frac{4s}{N}$ is the mass-critical exponent, and $q = \frac{p}{2}$. The parameter $\mu\in\mathbb{R}$ is the chemical potential, which enters as the Lagrange multiplier associated with the mass constraint, whereas $a_i>0\ (i=1,2)$ and $\beta>0$ measure the strengths of the intraspecies and interspecies interactions, respectively; the positivity of these coefficients means that both types of interaction are attractive.

	System \eqref{eq:sys} arises from a fundamental quantum mechanical model describing multicomponent Bose-Einstein condensates (BECs). Historically, the mathematical description of BECs, nonlinear optics, and many-body quantum physics has heavily relied on the classical (local) nonlinear Schr\"{o}dinger and Gross-Pitaevskii (GP) systems \cite{Dalfovo1999}. In particular, two-component GP systems model binary BECs with intraspecies and interspecies interactions, where the sign and strength of the interactions determine whether the condensate is stable, collapses, or undergoes phase separation \cite{Ho1996}. A natural nonlocal generalization is obtained when the particles in the condensate perform L\'{e}vy flights rather than Brownian motion: then the underlying kinematics are governed by the fractional Laplacian $(-\Delta)^s$ \cite{Laskin2000, Laskin2002}, which leads precisely to system \eqref{eq:sys}. In the study of such physical models, seeking \textit{normalized solutions}, namely solutions with a prescribed $L^2$ mass, has become a major theme because the mass directly corresponds to the number of particles, while the chemical potential naturally appears as a Lagrange multiplier. Before proceeding to the coupled fractional setting, it is instructive to recall the well-understood single-component theory, which provides the critical mass \(a^*\) and the ground state \(Q\) as basic building blocks.
	
	The theory for a single component is by now well understood. For the GP equation in $\mathbb R^2$ with attractive interaction, Guo and Seiringer \cite{GuoSeiringer2014} studied the mass concentration as the interaction strength $a$ approaches the critical value $a^*=\|Q\|_2^2$, where $Q$ is the unique positive radial ground state of
	\begin{equation}
	-\Delta u+u=u^3\quad\text{in }\mathbb R^2.
	\end{equation}
	They proved that minimizers exist if and only if $0\le a<a^*$, and that, as $a\nearrow a^*$, the minimizers blow up and concentrate at a global minimum of the trapping potential. Their proof relies on the sharp Gagliardo--Nirenberg inequality and delicate energy estimates. Later, Guo, Wang, Zeng and Zhou \cite{GuoWangZengZhou2018} extended this analysis to multi-well potentials, showing that ground states concentrate and blow up at the center of the incircle of the well with the largest inradius, and they established symmetry breaking. For ring-shaped potentials, Guo, Zeng and Zhou \cite{GuoZengZhou2016} proved that symmetry breaking occurs: for small $a$ the minimizer is unique and radially symmetric, whereas for $a$ close to $a^*$ there are infinitely many non-radial minimizers concentrating at points on the circular bottom of the potential well. These works rely on the sharp Gagliardo--Nirenberg inequality in the local case, refined energy estimates, blow-up analysis, and the moving plane method.
	
	These results have been extended to fractional Schr\"odinger equations. In the fractional setting, the corresponding ground state $Q$ of
	\begin{equation}
	(-\Delta)^su+\frac{2s}{N}u=|u|^{\frac{4s}{N}}u\quad\text{in }\mathbb R^N
	\end{equation}
	decays only polynomially, $Q(x)=O(|x|^{-N-2s})$, rather than exponentially, and the optimal constant in the fractional Gagliardo--Nirenberg inequality is more delicate. The uniqueness of the positive radial ground state $Q$ was proved by Frank and Lenzmann \cite{FrankLenzmann2013} in one dimension and by Frank, Lenzmann and Silvestre \cite{FrankLenzmannSilvestre2016} in higher dimensions, using nonlocal moving plane methods and the extension problem of Caffarelli and Silvestre. For fractional Schr\"odinger equations with trapping potentials, Du, Tian, Wang and Zhang \cite{DuTianWangZhang2019} considered
	\begin{equation}
	(-\Delta)^s u+V(x)u=\mu u+a|u|^{p-2}u,\quad x\in\mathbb R^N,
	\end{equation}
	with $0<s<1$, $N\ge2$, and $p=2+4s/N$. They established the existence and non-existence of $L^2$-normalized solutions by using a constrained variational method. They proved that minimizers exist if and only if $0\le a<a^*$, where $a^*=\|Q\|_2^{4s/N}$, and they analyzed the concentration behavior as $a\nearrow a^*$. The key tools are the fractional Gagliardo--Nirenberg inequality and compactness results in weighted Sobolev spaces. For fractional NLS with combined nonlinearities, Luo and Zhang \cite{LuoZhang2020} studied
	\begin{equation}
	(-\Delta)^s u=\lambda u+\mu|u|^{q-2}u+|u|^{p-2}u\quad\text{in }\mathbb R^N,
	\end{equation}
	with prescribed mass $\int|u|^2=a^2$, and obtained existence and non-existence results by means of fiber maps, Pohozaev identities, and minimax arguments. For the mass-supercritical case, Peng and Xia \cite{PengXia2021} studied
	\begin{equation}
	(-\Delta)^s u+V(x)u+\omega u=|u|^{p-2}u\quad\text{in }\mathbb R^N,
	\end{equation}
	with $p\in(2+4s/N,2_s^*)$, and proved the existence of normalized solutions under suitable assumptions on the potential, using a new minimax argument and a splitting lemma for the nonlocal operator.
	
	For systems, the situation is considerably richer. Two natural mass constraints appear in the literature. The first is the \emph{separate constraint}
	\begin{equation}
	\int_{\mathbb R^N}|u_1|^2\,dx=\int_{\mathbb R^N}|u_2|^2\,dx=1,
	\end{equation}
	under which each component has a prescribed mass. The second is the \emph{total-mass constraint}
	\begin{equation}\label{sum}
	\int_{\mathbb R^N}\bigl(|u_1|^2+|u_2|^2\bigr)\,dx=1,
	\end{equation}
	under which only the sum of the masses is fixed, and mass may be transferred between the two components. These two constraints lead to substantially different variational structures, critical thresholds, and blow-up phenomena.
	
	For two-component systems with separate mass constraints, Bartsch, Jeanjean and Soave \cite{BartschJeanjeanSoave2016} studied the coupled cubic Schr\"odinger system in $\mathbb R^3$,
	\begin{equation}
	\begin{cases}
		-\Delta u-\lambda_1 u=\mu_1 u^3+\beta u v^2,\\
		-\Delta v-\lambda_2 v=\mu_2 v^3+\beta u^2 v,
	\end{cases}
	\end{equation}
	with $\int_{\mathbb R^3}u^2=a_1^2$, $\int_{\mathbb R^3}v^2=a_2^2$. They proved the existence of positive normalized solutions for different ranges of the coupling parameter $\beta>0$, using a two-dimensional linking argument, a minimax principle, and a Pohozaev-type constraint. In $\mathbb R^2$, Guo, Zeng and Zhou \cite{GuoZengZhou2017} studied the system
	\begin{equation}
	\begin{cases}
		-\Delta u_1+V_1(x)u_1=\mu_1 u_1+b_1u_1^3+\beta u_2^2u_1,\\
		-\Delta u_2+V_2(x)u_2=\mu_2 u_2+b_2u_2^3+\beta u_1^2u_2,
	\end{cases}
	\end{equation}
	with $\int_{\mathbb R^2}u_1^2=\int_{\mathbb R^2}u_2^2=1$. They obtained a partial classification of the existence and non-existence of minimizers: if $0<b_1,b_2<a^*$ and $\beta<\sqrt{(a^*-b_1)(a^*-b_2)}$, minimizers exist; if $b_1>a^*$ or $b_2>a^*$ or $\beta>(a^*-b_1)/2+(a^*-b_2)/2$, no minimizer exists; the complete classification in the intermediate region $\sqrt{(a^*-b_1)(a^*-b_2)}\le\beta\le (a^*-b_1)/2+(a^*-b_2)/2$ remains open. They also proved that, as $(b_1,b_2)\nearrow (a^*-\beta,a^*-\beta)$, each component blows up at the same point and converges to the same profile after scaling, and they obtained the optimal blow-up rate. The proof uses an auxiliary minimization problem, refined energy estimates, and a blow-up analysis. In a subsequent paper \cite{GuoZengZhou2018}, they considered the case of attractive intraspecies and repulsive interspecies interactions ($\beta<0$) and obtained the blow-up behavior of ground states, showing that each component concentrates at a global minimum of the associated trapping potential.
	
	When the constraint is the total mass, the situation changes significantly because mass can be transferred between the two components. Guo, Li, Wei and Zeng \cite{GuoLiWeiZeng2019} studied the two-component attractive BEC system
	\begin{equation}
	\begin{cases}
		-\Delta u_1+V_1(x)u_1=\mu u_1+a_1u_1^3+\beta u_2^2u_1,\\
		-\Delta u_2+V_2(x)u_2=\mu u_2+a_2u_2^3+\beta u_1^2u_2,
	\end{cases}
	\end{equation}
	in $\mathbb R^2$ with $\int_{\mathbb R^2}(|u_1|^2+|u_2|^2)\,dx=1$. They gave a complete classification of the existence and non-existence of ground states, and analyzed uniqueness and symmetry breaking as $\beta\nearrow \beta^*=a^*+\sqrt{(a^*-a_1)(a^*-a_2)}$. Their methods include a refined Gagliardo--Nirenberg inequality, energy estimates, and a local uniqueness argument based on the implicit function theorem.
	
	For fractional systems with separate constraints, Liu, Zhang and Zuo \cite{LiuZhangZuoCPAA} studied the one-dimensional half-Laplacian system
	\begin{equation}
	\begin{cases}
		\sqrt{-\Delta}u_1+V_1(x)u_1=\mu_1u_1+a_1u_1^3+\beta u_2^2u_1,\\
		\sqrt{-\Delta}u_2+V_2(x)u_2=\mu_2u_2+a_2u_2^3+\beta u_1^2u_2,
	\end{cases}
	\quad x\in\mathbb R,
	\end{equation}
	with $\int_{\mathbb R}u_1^2=\int_{\mathbb R}u_2^2=1$. They established a criterion for the existence and non-existence of minimizers by solving an auxiliary minimization problem, and they derived delicate estimates for the energy functional. Their results show that if $0<a_1,a_2<a^*$ and $\beta<\beta_-:=\sqrt{(a^*-a_1)(a^*-a_2)}$, a minimizer exists, whereas if $a_1>a^*$ or $a_2>a^*$ or $\beta>\beta_+:=(a^*-a_1)/2+(a^*-a_2)/2$, no minimizer exists. The proof uses a refined fractional Gagliardo--Nirenberg inequality and a concentration-compactness argument. In a subsequent work \cite{LiuZhangZuoUniqueness}, they studied the uniqueness and concentration behavior of solutions for the same system under the separate constraint. They proved the uniqueness of nonnegative minimizers for small $|(a_1,a_2,\beta)|$ via the implicit function theorem, and analyzed the blow-up behavior as $(a_1,a_2)\nearrow(a^*-\beta,a^*-\beta)$. They showed that each component blows up and concentrates at a flattest common minimum point of $V_1$ and $V_2$, and obtained the optimal blow-up rate
	\begin{equation}
		\varepsilon_k=\Bigl(a^*-\frac{a_{1k}+a_{2k}+2\beta}{2}\Bigr)^{\frac1{p_0+1}},
	\end{equation}
	where $p_0$ is determined by the order of the flattest common minimum. The proof combines energy estimates, Pohozaev identities, and nonlocal regularity theory.
	
	It is important to observe that the classification obtained in \cite{LiuZhangZuoCPAA} for the separate constraint is not complete. Indeed, while the existence and non-existence of minimizers are settled when $\beta<\beta_-$ and when $\beta>\beta_+$, the intermediate region $\beta_-\le\beta\le\beta_+$ is only partially understood: for $s=\frac 12$, Theorems 1.5 and 1.6 of \cite{LiuZhangZuoCPAA} provide some partial results under additional assumptions, whereas for $s\ne\frac 12$ the complete classification remains an open problem (see Remark 1.3 therein). Moreover, the blow-up analysis in \cite{LiuZhangZuoUniqueness} is restricted to the separate constraint, and the optimal rate is obtained under specific assumptions on the trapping potentials.
	
To the best of our knowledge, however, neither a complete classification nor a description of the concentration behavior of ground states of the coupled fractional system \eqref{eq:sys} under the total-mass constraint \eqref{sum} is available in the literature. Inspired by the above works, this paper is devoted to two issues: the complete classification of ground states of \eqref{eq:sys} for arbitrary $s\in(0,1)$ and $N\ge1$, and the concentration behavior of the minimizers in the one-dimensional half-Laplacian case $N=1$, $s=\frac 12$, for trapping potentials $V_i(x)$ $(i=1,2)$ satisfying

\begin{itemize}
	\item[$(\Omega)$] $0 \le V_i(x) \in L_{loc}^\infty(\mathbb{R}^N)$, $\lim_{|x|\to\infty} V_i(x) = \infty$, and $\inf_{x \in \mathbb{R}^N} \big(V_1(x) + V_2(x)\big) = 0$ is attained.
\end{itemize}

We begin with some notation. Let the Sobolev space
\begin{equation}
	H^s(\mathbb R^N):=\left\{u\in L^2(\mathbb R^N):\|u\|^2_{H^{s}(\mathbb{R}^N)} := \int_{\mathbb R^N} \big(1+|\xi|^2\big)^s |\mathcal F u(\xi)|^2 d\xi <\infty\right\}.
\end{equation}
Here $\mathcal{F}u$ denotes the Fourier transform of $u$. Equipped with the inner product
\begin{equation}
	(u, v)_{H^s(\mathbb{R}^N)} := \int_{\mathbb{R}^N} \big(1+|\xi|^2\big)^s \overline{\mathcal{F}u(\xi)} \mathcal{F}v(\xi) d\xi,
\end{equation}
the space $H^s(\mathbb{R}^N)$ forms a Hilbert space. Moreover, the norm $\|\cdot\|_{H^s(\mathbb{R}^N)}$ is equivalent to the one given by
\begin{equation}
	\|u\|_{H^s(\mathbb{R}^N)}^2 = \int_{\mathbb{R}^N} \left[ u (-\Delta)^s u + u^2 \right] dx = \int_{\mathbb{R}^N} \left[ \left| (-\Delta)^{\frac{s}{2}} u \right|^2 + u^2 \right] dx.
\end{equation}

	Ground states of system \eqref{eq:sys} can be studied equivalently through the minimizers of the following constrained variational problem (see Proposition \ref{Proposition A.1.}):
\begin{equation}\label{eq:min_problem}
	e(a_1, a_2, \beta) := \inf_{(u_1, u_2) \in \mathcal{M}} E_{a_1, a_2, \beta}(u_1, u_2),
\end{equation}
where the associated Gross-Pitaevskii (GP) energy functional is given by
\begin{equation}\label{eq:energy}
	\begin{aligned}
		E_{a_1, a_2, \beta}(u_1, u_2)
		= &\sum_{i=1}^2 \int_{\mathbb{R}^N} \left( |(-\Delta)^{\frac s2}u_i|^2 + V_i(x)u_i^2 \right) dx \corr{-} \frac{2}{p}\sum_{i=1}^2 \int_{\mathbb{R}^N} a_i|u_i|^p dx \\
		& - \frac{4}{p} \beta \int_{\mathbb{R}^N} |u_1|^q |u_2|^q dx,\quad (u_1,u_2)\in \mathcal H_1\times \mathcal H_2.\\
	\end{aligned}
\end{equation}
Here, the mass constraint $\mathcal{M}$ is defined by
\begin{equation}\label{mass}
	\mathcal{M}:=\left\{(u_1,u_2)\in \mathcal{H}_1\times \mathcal{H}_2: \int_{\mathbb R^N} \big(u_1^2+u_2^2\big)dx=1\right\}
\end{equation}
with
\begin{equation}\label{eq:1-11}
	\mathcal{H}_i=\left\{u\in H^{s}(\mathbb R^N): \int_{\mathbb R^N} V_i(x)u^2 dx<\infty\right\},\quad i=1,2,
\end{equation}
and
\begin{equation}\label{eq:1-12}
	\|u\|_{\mathcal{H}_i}^2=\int_{\mathbb R^N} \left[\big|(-\Delta)^{\frac s2}u\big|^2 + V_i(x)u^2\right]dx.
\end{equation}

	In order to investigate the coupled minimization problem \eqref{eq:min_problem}, it is essential to recall several fundamental results concerning the single-component fractional equation. The classical scalar field equation associated with system \eqref{eq:sys} in the mass-critical case is known as
\begin{equation}\label{singal}
	(-\Delta)^s u + \frac{2s}{N} u - |u|^{\frac{4s}{N}}u = 0 \quad \text{in } \mathbb{R}^N, \quad u \in H^s(\mathbb{R}^N).
\end{equation}
As is well known (see \cite{FrankLenzmann2013,FrankLenzmannSilvestre2016}), there exists a unique (up to translations) positive radial ground state $Q\in\Hd$ of equation \eqref{singal}. Moreover, $Q$ has the decay estimate
\begin{equation}\label{decay}
	Q(x)=O\bigl( |x|^{-N-2s}\bigr) \quad \text{as} \quad |x|\to\infty.
\end{equation}
The following Gagliardo-Nirenberg inequality holds with optimal constant:
\begin{equation}
	\intR \abs{u}^p dx \leq \frac{N+2s}{N\norm{Q}{2}^{\frac{4s}{N}}} \norm{(-\Delta)^{\frac s2} u}{2}^2 \norm{u}{2}^{\frac{4s}{N}} \quad \forall u\in\Hd.
\end{equation}
Define $a^* := \norm{Q}{2}^{\frac{4s}{N}}$. This constant plays the role of a critical mass for the single equation. Thus, we have
\begin{equation}\label{lem:single_gns}
	\intR \abs{u}^p dx \leq \frac{p}{2a^*} \norm{(-\Delta)^{\frac s2} u}{2}^2 \norm{u}{2}^{\frac{4s}{N}} \quad \forall u\in\Hd.
\end{equation}

	We now state our main results. The first result concerns the complete classification of existence and non-existence of the minimizers for \eqref{eq:min_problem}. To deal with the total-mass constraint, the main mathematical challenge lies in providing a lower bound for the energy functional. Instead of artificially constructing a new coupled Gagliardo-Nirenberg inequality, we introduce an auxiliary intensity function $\Lambda(a_1,a_2,\beta)$ and characterize the exact critical coupling threshold $\beta^*$ as follows:
\begin{equation}\label{eq:Lambda}
	\Lambda(a_1, a_2, \beta) := \max_{\theta \in [0,1]} \left( a_1\theta^q + a_2(1-\theta)^q + 2\beta\theta^{q/2}(1-\theta)^{q/2} \right),
\end{equation}
\begin{equation}\label{eq:beta}
	\beta^* := \min_{\theta \in (0,1)} \frac{a^* - a_1\theta^q - a_2(1-\theta)^q}{2\theta^{q/2}(1-\theta)^{q/2}}.
\end{equation}

\begin{remark}\label{state}
	We note that for $0 \leq a_1, a_2 < a^*$, the critical coupling threshold $\beta^*$ given by \eqref{eq:beta} is well-defined. Indeed, the function to be minimized tends to $+\infty$ as $\theta \to 0^+$ and $\theta \to 1^-$, ensuring that the minimum is attained in the interior $(0,1)$. Furthermore, in this parameter regime, it is straightforward to verify the algebraic equivalence that $\Lambda(a_1, a_2, \beta) < a^*$ if and only if $\beta < \beta^*$.
\end{remark}

Our first main result establishes the existence and non-existence of the minimizers of the problem \eqref{eq:min_problem} for arbitrary $0<s<1$ and $N \ge 1 $.

\begin{theorem} \label{thm:main}
	Suppose that $(\Omega)$ holds. Then the following statements hold.\\
	(i) If $0 \leq a_1, a_2 < a^*$ and $0 < \beta < \beta^*$, then problem \eqref{eq:min_problem} has at least one minimizer. Moreover, it holds
	\begin{equation} \label{eq:energy_limit}
		\lim_{\beta \to \beta^*} e(a_1,a_2,\beta) = e(a_1,a_2,\beta^*) = 0, \quad 0 \leq a_1, a_2 < a^*.
	\end{equation}
	(ii) If $a_1 \geq a^*$ or $a_2 \geq a^*$ or $\beta \geq \beta^*$, then problem \eqref{eq:min_problem} has no minimizer.
\end{theorem}
\begin{remark}
	We remark that the existence result also holds true for the system with repulsive interspecies interaction, i.e., $\beta<0$, and for the system with $\beta=0$; this can be seen from the proof of Theorem 1.1 (i). This complete classification identifies the exact parameter regimes in which the total-mass-constrained fractional BEC system remains stable. In particular, it covers the degenerate cases $a_1=0$ or $a_2=0$ and yields a complete classification for all $\beta\in\mathbb{R}$, in contrast with the partial classification available in the separate-constraint setting of \cite{LiuZhangZuoCPAA}.
\end{remark}

Our second goal is to study the blow-up behavior and the optimal blow-up rate for the one-dimensional half-Laplacian system under the total-mass constraint, obtaining results analogous to the last two results of \cite{GuoLiWeiZeng2019}. More precisely, we consider
\begin{equation}
\begin{cases}
	\sqrt{-\Delta}\,u_1+V_1(x)u_1=\mu u_1+a_1u_1^3+\beta u_2^2u_1,\\[1mm]
	\sqrt{-\Delta}\,u_2+V_2(x)u_2=\mu u_2+a_2u_2^3+\beta u_1^2u_2,
\end{cases}
\qquad x\in\mathbb R,
\end{equation}
with $\int_{\mathbb R}(|u_1|^2+|u_2|^2)\,dx=1$, which implies $p=4$ and $q=2$. By \eqref{eq:beta}, one can easily verify that the minimum is attained at the optimal mass fraction
\begin{equation}\label{eq:3-1}
	\gamma:=\frac{\sqrt{a^*-a_2}}{\sqrt{a^*-a_1}+\sqrt{a^*-a_2}},
\end{equation}
which directly yields $\beta^*=a^*+\sqrt{(a^*-a_1)(a^*-a_2)}$. The mass fraction $\gamma$ and the threshold $\beta^*$ play a crucial role in the proofs of the following theorems.

\begin{theorem}\label{con}
	Suppose $(\Omega)$ holds and let $(u_{1\beta}, u_{2\beta})$ be a nonnegative minimizer of $e(a_1, a_2, \beta)$ with $N=1$ and $s=\frac 12$, where $0 < a_1, a_2 < a^*$ and $0 < \beta < \beta^*$. Then for any sequence $\{\beta_k\}$ with $\beta_k \to \beta^*$ as $k \to \infty$, there exists a subsequence of $\{\beta_k\}$, still denoted by $\{\beta_k\}$, such that for $i=1,2$ and every $k$, $u_{i\beta_k}$ has at least one global maximum point $z_{i\beta_k}$, and
	\begin{equation}\label{eq:2-16}
		\lim_{k\to\infty} z_{1\beta_k} = \lim_{k\to\infty} z_{2\beta_k} = x_0 \quad \text{for some } x_0 \in \mathcal{V},
	\end{equation}
	where
	\begin{equation}
		\mathcal{V} := \bigl\{ x \in \mathbb{R} : V_1(x) = V_2(x) = 0 \bigr\}.
	\end{equation}
	Moreover, we have
	\begin{equation}\label{eq:2-18}
		\begin{cases}
			\displaystyle\lim_{k\to\infty} \varepsilon_{\beta_k}^{\frac{1}{2}} u_{1\beta_k}\bigl(\varepsilon_{\beta_k} x + z_{1\beta_k}\bigr) = \frac{\sqrt{\gamma}}{\|Q\|_2} Q(x), \\[1em]
			\displaystyle\lim_{k\to\infty} \varepsilon_{\beta_k}^{\frac{1}{2}} u_{2\beta_k}\bigl(\varepsilon_{\beta_k} x + z_{2\beta_k}\bigr) = \frac{\sqrt{1-\gamma}}{\|Q\|_2} Q(x),
		\end{cases}
	\end{equation}
	strongly in $H^{\frac{1}{2}}(\mathbb{R})$, where $\gamma$ is given by \eqref{eq:3-1} and $Q$ is the positive radially symmetric ground state of \eqref{singal}; here $\varepsilon_{\beta_k} > 0$ satisfies
	\begin{equation}\label{epsilo}
		\varepsilon_{\beta_k} := \biggl( \int_{\mathbb{R}} \bigl( |(-\Delta)^{\frac{1}{4}} u_{1\beta_k}|^2 + |(-\Delta)^{\frac{1}{4}} u_{2\beta_k}|^2 \bigr) dx \biggr)^{-1} \to 0 \quad \text{as } k \to \infty,
	\end{equation}
	\begin{equation}\label{eq:2-20}
		\lim_{k\to\infty} \frac{|z_{1\beta_k} - z_{2\beta_k}|}{\varepsilon_{\beta_k}} = 0,
	\end{equation}
	and
	\begin{equation}\label{eq:2-21}
		\lim_{k\to\infty} \frac{\beta^* - \beta_k}{\varepsilon_{\beta_k}} = 0.
	\end{equation}
\end{theorem}

\begin{remark}
	Unlike \cite{GuoLiWeiZeng2019}, we do not establish the uniqueness of the global maximum point of each component; we only prove that at least one global maximum point exists. The uniqueness remains an open problem, mainly due to the lack of sufficient elliptic regularity for the fractional operator, which prevents us from proving the smoothness of the minimizer sequence (see also Remark 1.4 in \cite{LiuZhangZuoUniqueness}).
\end{remark}

To carry out a finer analysis of the concentration behavior of nonnegative minimizers, we henceforth assume that the trapping potentials $V_i(x)$ have the form
\begin{equation}\label{eq:potential_form}
	V_i(x) = h_i(x) \prod_{j=1}^{m_i} |x - x_{ij}|^{q_{ij}}, \quad q_{ij} > 0 \text{ and } m_i \in \mathbb{N}^+, \ i = 1, 2,
\end{equation}
where $h_i(x) \in C_{loc}^\alpha(\mathbb{R})$ satisfies $\frac{1}{C} < h_i(x) < C$ for some positive constant $C$ and $x_{ij} \neq x_{ik}$ for $j \neq k$. We also assume, without loss of generality, that there is an integer $b$ with $1 \le b \le \min\{m_1, m_2\}$ such that
\begin{equation}\label{eq:common_zeros_condition}
	\begin{cases}
		x_{1j} = x_{2j} := x_j, & \text{where } j = 1, \dots, b; \\
		x_{1j} \neq x_{2k}, & \text{where } j \in \{b+1, \dots, m_1\} \text{ and } k \in \{b+1, \dots, m_2\}.
	\end{cases}
\end{equation}
Notice that \eqref{eq:common_zeros_condition} implies
\begin{equation}\label{eq:set_V}
	\mathcal{V} = \{ x \in \mathbb{R} : V_1(x) = V_2(x) = 0 \} = \{x_1, \dots, x_b\}.
\end{equation}

Additionally, we define $\lambda_j$ for $j = 1, \dots, b$ by
\begin{equation}\label{eq:lambda_j}
	\lambda_j := \begin{cases}
		\displaystyle\lim_{x \to x_j} \frac{\gamma V_1(x)}{|x - x_j|^{q_{1j}}} \int_{\mathbb{R}} |y|^{q_{1j}} Q^2(y) \, dy & \text{if } q_{1j} < q_{2j}, \\[12pt]
		\displaystyle\lim_{x \to x_j} \frac{\gamma V_1(x) + (1-\gamma)V_2(x)}{|x - x_j|^{q_{1j}}} \int_{\mathbb{R}} |y|^{q_{1j}} Q^2(y) \, dy & \text{if } q_{1j} = q_{2j}, \\[12pt]
		\displaystyle\lim_{x \to x_j} \frac{(1-\gamma) V_2(x)}{|x - x_j|^{q_{2j}}} \int_{\mathbb{R}} |y|^{q_{2j}} Q^2(y) \, dy & \text{if } q_{1j} > q_{2j},
	\end{cases}
\end{equation}
where $\gamma$ is the optimal mass fraction given by \eqref{eq:3-1}. Define
\begin{equation}\label{eq:q_j_q_0}
	q_j := \min\{q_{1j}, q_{2j}\}, \quad j = 1, \dots, b; \quad q_0 := \max_{1 \le j \le b} q_j,
\end{equation}
so that
\begin{equation}\label{eq:set_Z}
	Z := \{ x_j : q_j = q_0, j = 1, \dots, b \} \subset \mathcal{V}.
\end{equation}
Set
\begin{equation}\label{eq:lambda_0}
	\lambda_0 := \min_{j \in \Gamma} \lambda_j, \quad \text{where } \Gamma := \{j : x_j \in Z\},
\end{equation}
and let
\begin{equation}\label{eq:set_Z_0}
	Z_0 := \{ x_j \in Z : \lambda_j = \lambda_0 \}
\end{equation}
denote the set of the flattest common minimum points of $V_i(x)$ for $i = 1, 2$. Our third main result concerns the precise concentration behavior of nonnegative minimizers for \eqref{eq:min_problem} with a polynomial trapping potential defined above.

\begin{theorem}\label{thm:1.3}
	Suppose that $V_i(x)$ ($i = 1, 2$) satisfy \eqref{eq:potential_form} and \eqref{eq:common_zeros_condition}, and let $(u_{1\beta}, u_{2\beta})$ be a nonnegative minimizer of $e(a_1, a_2, \beta)$ with $N=1$ and $s=\frac 12$, where $0 < a_1, a_2 < a^*$ and $0 < \beta < \beta^*$. Moreover, if $0<q_0<1$, then for any sequence $\{\beta_k\}$ with $\beta_k \to \beta^*$ as $k \to \infty$, there exists a subsequence of $\{\beta_k\}$, still denoted by $\{\beta_k\}$, such that for $i=1,2$ and every $k$, $u_{i\beta_k}$ has at least one global maximum point $z_{i\beta_k}$, and
	\begin{equation}\label{eq:z_limit}
		\lim_{k\to\infty} \frac{z_{1\beta_k} - x_{j_0}}{\epsilon_{\beta_k}} = \lim_{k\to\infty} \frac{z_{2\beta_k} - x_{j_0}}{\epsilon_{\beta_k}} = 0 \quad \text{for some } x_{j_0} \in Z_0,
	\end{equation}
	where $\epsilon_{\beta_k} > 0$ is given by
	\begin{equation}\label{eq:epsilon_beta_k}
		\epsilon_{\beta_k} := \left[ \frac{2\gamma(1-\gamma)}{q_0\lambda_0} (\beta^* - \beta_k) \right]^{\frac{1}{q_0+1}}.
	\end{equation}
	Moreover, we have
	\begin{equation}\label{eq:profile_limit}
		\begin{cases}
			\lim_{k\to\infty} \epsilon_{\beta_k}^{\frac{1}{2}} u_{1\beta_k}(\epsilon_{\beta_k} x + z_{1\beta_k}) = \frac{\sqrt{\gamma}}{\|Q\|_2} Q(x), \\
			
			\lim_{k\to\infty} \epsilon_{\beta_k}^{\frac{1}{2}} u_{2\beta_k}(\epsilon_{\beta_k} x + z_{2\beta_k}) = \frac{\sqrt{1-\gamma}}{\|Q\|_2} Q(x),
		\end{cases}
	\end{equation}
	strongly in $H^{\frac{1}{2}}(\mathbb{R})$ and
	\begin{equation}\label{eq:energy_limit_1}
		\lim_{k\to\infty} \epsilon_{\beta_k} \int_{\mathbb{R}} \left( |(-\Delta)^{\frac{1}{4}} u_{1\beta_k}|^2 + |(-\Delta)^{\frac{1}{4}} u_{2\beta_k}|^2 \right) dx = 1.
	\end{equation}
\end{theorem}

Theorem \ref{thm:1.3} shows that, despite the long-range algebraic interaction and the strong attractive coupling, both components retain the mass ratio $\gamma:(1-\gamma)$ and concentrate simultaneously at the flattest common minimum.

The main difficulties of this work are three aspects. First, the passage from the classical Laplacian to the fractional operator $(-\Delta)^s$ with arbitrary $s\in(0,1)$ and dimension $N\ge1$ is not merely technical. The nonlocal nature of $(-\Delta)^s$ destroys the usual maximum principle and ODE techniques; one has to rely on nonlocal elliptic regularity, De Giorgi--Nash--Moser theory, and the positivity of the resolvent kernel. Moreover, the critical exponents $p=2+4s/N$ and $q=1+2s/N$ depend on both $s$ and $N$, so the structure of the nonlinearities changes with the dimension, and the ground state $Q$ decays polynomially rather than exponentially, which makes the standard test-function estimates and the derivation of optimal energy bounds much more delicate. Second, the change from the separate constraint to the total-mass constraint is substantial: mass may be transferred between the two components, and the energy functional no longer admits an obvious lower bound; this difficulty is overcome by means of the auxiliary intensity function $\Lambda(a_1,a_2,\beta)$, which yields the exact critical coupling threshold $\beta^*$ (defined in \eqref{eq:beta}) for arbitrary $s$ and $N$. It is only in the specific one-dimensional half-Laplacian case ($N=1$, $s=1/2$) that this threshold simplifies to $\beta^* = a^* + \sqrt{(a^*-a_1)(a^*-a_2)}$. This auxiliary method naturally covers the degenerate cases ($a_1=0$ or $a_2=0$) and the repulsive interspecies regime ($\beta<0$), providing a complete classification for all $\beta\in\mathbb{R}$. Finally, the blow-up analysis requires careful treatment of the mutual interaction to show that both components concentrate simultaneously at the flattest common minimum.

Our analysis combines several ingredients: refined fractional Gagliardo--Nirenberg inequalities, constrained minimization, concentration-compactness arguments, Pohozaev identities, energy estimates, and nonlocal regularity theory. The paper is organized as follows. In Section 2, we establish the existence and non-existence criteria (Theorem \ref{thm:main}). Section 3 is devoted to the proof of Theorem \ref{con} on the general limit behavior, while Section 4 finishes the proof of Theorem \ref{thm:1.3} on the explicit blow-up profiles under polynomial trapping potentials.

\section{Existence and non-existence results}

This section is devoted to establishing the existence and non-existence of minimizers for problem \eqref{eq:min_problem}. We begin by recalling several preliminary results essential to our subsequent analysis.

\begin{lemma} \label{lem:compact_embedding}
Under condition $(\Omega)$, the embedding $\mathcal{H}_1\times\mathcal{H}_2\hookrightarrow L^r(\R^N)\times L^r(\R^N)$ is compact for any $2\leq r<2_s^*$, where $2_s^*=\frac{2N}{N-2s}$ for $N>2s$ and $2_s^*=\infty$ for $N\leq2s$.
\end{lemma}

\begin{proof}
The proof is standard; see e.g. \cite[Lemma 3.1]{DuTianWangZhang2019}. The condition $V_i(x)\to\infty$ ensures that bounded sequences in $\mathcal{H}_i$ are tight, and local compactness of $H^s$ gives strong convergence in $L^r$ after extracting a subsequence.
\end{proof}
\begin{remark}
	Note that the exponent $p=2+\frac{4s}{N}$ satisfies $2\leq p<2_s^*$, hence $p$ lies in the compact range.
\end{remark}
\begin{lemma} \label{holder}
Let $u_1, u_2 \in L^p(\mathbb{R}^N)$. Then,
\begin{equation} \label{eq:holder_coupled}
\int_{\mathbb{R}^N} |u_1|^q |u_2|^q dx \leq \left( \int_{\mathbb{R}^N} |u_1|^p dx \int_{\mathbb{R}^N} |u_2|^p dx \right)^{1/2},
\end{equation}
and the equality holds if and only if $|u_1|^p = \kappa |u_2|^p$ for some $\kappa > 0$.
\end{lemma}

\begin{lemma} \label{lem:convergence}
If $(u_{1n}, u_{2n}) \rightharpoonup (u_{10}, u_{20})$ weakly in $\mathcal{H}_1\times\mathcal{H}_2$, then, by Lemma \ref{lem:compact_embedding}, $u_{in}\to u_{i0}$ strongly in $L^r(\mathbb{R}^N)$ for every $2\le r<2_s^*$, and hence
\begin{equation}
\lim_{n \to \infty} \int_{\mathbb{R}^N} |u_{1n}|^q |u_{2n}|^q dx = \int_{\mathbb{R}^N} |u_{10}|^q |u_{20}|^q dx,
\end{equation}
\begin{equation}
\lim_{n \to \infty} \int_{\mathbb{R}^N} |u_{in}|^p dx = \int_{\mathbb{R}^N} |u_{i0}|^p dx, \quad i=1,2.
\end{equation}
\end{lemma}

\begin{proof}
Since $u_{in}\to u_{i0}$ strongly in $L^p(\mathbb{R}^N)$, the sequences $\{u_{in}\}$ are uniformly bounded in $L^p(\mathbb{R}^N)$. Note that $p/q=2$, which allows us to work in $L^2$. Using the elementary inequality
\begin{equation}
\bigl||a|^q-|b|^q\bigr| \le C_q\bigl(|a|^{q-1}+|b|^{q-1}\bigr)|a-b|,\qquad \forall a,b\in\mathbb{R},
\end{equation}
and Hölder's inequality, we obtain
\begin{equation}
\bigl\||u_{1n}|^q-|u_{10}|^q\bigr\|_{2}
\le C_q\bigl(\|u_{1n}\|_{p}^{q-1}+\|u_{10}\|_{p}^{q-1}\bigr)\|u_{1n}-u_{10}\|_{p}\to0 \quad \text{as} \ n \to \infty,
\end{equation}
and similarly $\bigl\||u_{2n}|^q-|u_{20}|^q\bigr\|_{2}\to0 $ as $n\to \infty $.

Now estimate the difference of the integrals, by Hölder's inequality, we have

\begin{equation}\label{nonlinear_term}
  \begin{aligned}
&\left|\int_{\mathbb{R}^N} |u_{1n}|^q|u_{2n}|^q\,dx-\int_{\mathbb{R}^N} |u_{10}|^q|u_{20}|^q\,dx\right|\\
&\le \int_{\mathbb{R}^N} \bigl||u_{1n}|^q-|u_{10}|^q\bigr|\,|u_{2n}|^q\,dx
   +\int_{\mathbb{R}^N} |u_{10}|^q\,\bigl||u_{2n}|^q-|u_{20}|^q\bigr|\,dx\\
&\le \bigl\||u_{1n}|^q-|u_{10}|^q\bigr\|_{2}\,\|u_{2n}\|_{p}^{q}
+ \bigl\||u_{2n}|^q-|u_{20}|^q\bigr\|_{2}\,\|u_{10}\|_{p}^{q}.
  \end{aligned}
\end{equation}
The proof is finished by letting $n\to \infty$ in \eqref{nonlinear_term}.
\end{proof}

\begin{lemma} \label{lem:refined_gns}
For any $(u_1,u_2)\in\Hd\times\Hd$, the following refined Gagliardo-Nirenberg inequality holds:
\begin{equation}\label{refine}
\begin{aligned}
\intR \left( \abs{u_1}^2 + \abs{u_2}^2 \right)^q dx &\leq \frac{p}{2a^*} \left( \norm{(-\Delta)^{\frac s2} u_1}{2}^2 + \norm{(-\Delta)^{\frac s2} u_2}{2}^2 \right)\\
&\quad \times \left( \intR (u_1^2+u_2^2)\,dx \right)^{\frac{2s}{N}}.
\end{aligned}
\end{equation}
\end{lemma}

\begin{proof}
	
Let $v=\sqrt{u_1^2+u_2^2}$. Then $v\geq0$ and $ \norm{v}{2}^2 = \intR (u_1^2+u_2^2) dx.$
Combining \eqref{lem:single_gns} and \eqref{convex}, we obtain

\begin{equation} \label{eq:refined_gns_final}
	\begin{aligned}
		\intR \left( \abs{u_1}^2 + \abs{u_2}^2 \right)^q dx
		&= \intR v^p dx
		\leq \frac{p}{2a^*} \norm{(-\Delta)^{\frac s2} v}{2}^2 \norm{v}{2}^{\frac{4s}{N}} \\
		&\leq \frac{p}{2a^*} \left( \norm{(-\Delta)^{\frac s2} u_1}{2}^2 + \norm{(-\Delta)^{\frac s2} u_2}{2}^2 \right) \\
&\quad \times \left( \intR (u_1^2+u_2^2)\,dx \right)^{\frac{2s}{N}}.
	\end{aligned}
\end{equation}
This completes the proof of \eqref{refine}.
\end{proof}

\begin{proofof}{Theorem \ref{thm:main}}
Our proof is divided into three parts.

\vspace{2mm}
(i) In the case $0 \leq a_1, a_2 < a^*$ and $0 < \beta < \beta^*$, we claim that \eqref{eq:min_problem} has at least one minimizer.
For any $(u_1, u_2) \in \mathcal{M}$, we introduce the total amplitude function $v(x) = \sqrt{u_1^2(x) + u_2^2(x)}$. On the support where $v(x) > 0$, we define the mass distribution fraction $\theta(x) = \frac{u_1^2(x)}{v^2(x)} \in [0, 1]$.
Because the exponents satisfy $p = 2q$, the nonlinear terms can be factorized pointwise:
\begin{equation}
\begin{aligned}
    \mathcal{N}(u_1, u_2) &:= a_1 |u_1|^p + a_2 |u_2|^p + 2\beta |u_1|^q |u_2|^q \\
    &= v(x)^p \left[ a_1 \theta(x)^q + a_2 (1-\theta(x))^q + 2\beta \theta(x)^{q/2}(1-\theta(x))^{q/2} \right].
    \end{aligned}
\end{equation}
By the definition of $\Lambda(a_1, a_2, \beta)$ in \eqref{eq:Lambda}, we have $\mathcal{N}(u_1, u_2) \le \Lambda(a_1, a_2, \beta) v(x)^p$ almost everywhere. Combining \eqref{mass}, \eqref{lem:single_gns} and \eqref{convex}, and recalling that $\|v\|_2=1$ on $\mathcal{M}$, we obtain
\begin{align}\label{eq:nonlinear}
    \int_{\mathbb{R}^N} \mathcal{N}(u_1, u_2) dx &\le \Lambda(a_1, a_2, \beta) \int_{\mathbb{R}^N} v^p dx \notag \\
    &\le \Lambda(a_1, a_2, \beta) \frac{p}{2a^*} \|(-\Delta)^{\frac s2} v\|_2^2 \|v\|_2^{\frac{4s}{N}} \notag \\
    &\le \frac{p \Lambda(a_1, a_2, \beta)}{2a^*} \sum_{i=1}^2 \|(-\Delta)^{\frac s2} u_i\|_2^2.
\end{align}
It follows from \eqref{eq:energy}, \eqref{eq:nonlinear} and Remark \ref{state}  that
\begin{equation}\label{eq:lower_bound}
\begin{aligned}
    E_{a_1, a_2, \beta}(u_1, u_2) &= \sum_{i=1}^2 \left( \|(-\Delta)^{\frac s2}u_i\|_2^2 + \int_{\mathbb{R}^N} V_i u_i^2 dx \right) - \frac{2}{p} \int_{\mathbb{R}^N} \mathcal{N}(u_1, u_2) dx \\
    &\ge \left( 1 - \frac{\Lambda(a_1, a_2, \beta)}{a^*} \right) \sum_{i=1}^2 \|(-\Delta)^{\frac s2} u_i\|_2^2 + \sum_{i=1}^2 \int_{\mathbb{R}^N} V_i u_i^2 dx\\
    &\geq0
\end{aligned}
\end{equation}
holds for any $(u_1,u_2)\in\mathcal{M}$. Consequently, $ E_{a_1, a_2, \beta}(u_1, u_2)$ is bounded from below in $\mathcal{M}$.
Let $\{(u_{1n}, u_{2n})\} \subset \mathcal{M}$ be a minimizing sequence for problem \eqref{eq:min_problem}, which satisfies
\begin{equation}\label{eq:min_seq}
	\int_{\mathbb{R}^N} (u_{1n}^2 + u_{2n}^2) dx = 1, \quad \lim_{n \to \infty} E_{a_1, a_2, \beta}(u_{1n}, u_{2n}) = e(a_1, a_2, \beta).
\end{equation}
Analogous to the derivation of \eqref{eq:lower_bound}, one can infer that
\begin{equation}
	E_{a_1, a_2, \beta}(u_{1n}, u_{2n}) \ge \left( 1 - \frac{\Lambda(a_1, a_2, \beta)}{a^*} \right) \sum_{i=1}^2 \|(-\Delta)^{\frac s2} u_{in}\|_2^2 + \sum_{i=1}^2 \int_{\mathbb{R}^N} V_i(x) u_{in}^2 dx.
\end{equation}
Passing to the limit $n \to \infty$ on both sides of the above inequality, we deduce from \eqref{eq:min_seq} that the minimizing sequence $\{(u_{1n}, u_{2n})\}$ is uniformly bounded in $\mathcal{H}_1 \times \mathcal{H}_2$. From Lemma \ref{lem:compact_embedding}, passing to a subsequence if necessary, there exists $(u_{10}, u_{20}) \in \mathcal{H}_1 \times \mathcal{H}_2$ such that
\begin{equation}
	\begin{cases}
		(u_{1n}, u_{2n}) \rightharpoonup (u_{10}, u_{20}) & \text{weakly in } \mathcal{H}_1 \times \mathcal{H}_2, \\
		(u_{1n}, u_{2n}) \to (u_{10}, u_{20}) & \text{strongly in } L^r(\mathbb{R}^N) \times L^r(\mathbb{R}^N), \quad 2 \le r < 2_s^*.
	\end{cases}
\end{equation}
Hence, we can deduce that
\begin{equation}
	\lim_{n \to \infty} \int_{\mathbb{R}^N} (u_{1n}^2 + u_{2n}^2) dx = 1 = \int_{\mathbb{R}^N} (u_{10}^2 + u_{20}^2) dx, \quad \text{i.e., } (u_{10}, u_{20}) \in \mathcal{M}.
\end{equation}
Furthermore, it follows from the weak lower semi-continuity of the norms $\|\cdot\|_{\mathcal{H}_i}$ and Lemma \ref{lem:convergence} that
\begin{equation}
	e(a_1, a_2, \beta) \le E_{a_1, a_2, \beta}(u_{10}, u_{20}) \le \lim_{n \to \infty} E_{a_1, a_2, \beta}(u_{1n}, u_{2n}) = e(a_1, a_2, \beta).
\end{equation}
This implies that $E_{a_1, a_2, \beta}(u_{10}, u_{20}) = e(a_1, a_2, \beta)$, establishing $(u_{10}, u_{20})$ as a minimizer of the problem \eqref{eq:min_problem}.

\vspace{2mm}
(ii) In the case $a_1 > a^*$ or $a_2 > a^*$ or $\beta > \beta^*$, we infer that \eqref{eq:min_problem} has no minimizers.
Let $\varphi \in C_0^\infty(\mathbb{R}^N)$ be a cut-off function satisfying
\begin{equation}\label{varphi}
\varphi(x) \equiv 1 \text{ for } |x| \leq 1, \quad \varphi(x) \equiv 0 \text{ for } |x| \geq 2, \quad 0 \leq \varphi(x) \leq 1 \text{ for all } x \in \mathbb{R}^N.
\end{equation}
Besides, for any $\sigma > 0$, set
\begin{equation} \label{eq:test_sequence}
\begin{cases}
\phi_{1\sigma}(x) := \frac{\sqrt{\theta} A_\sigma}{\|Q\|_2} \sigma^{N/2} \varphi(x-x_0) Q(\sigma(x-x_0)), \\
\phi_{2\sigma}(x) := \frac{\sqrt{1-\theta} A_\sigma}{\|Q\|_2} \sigma^{N/2} \varphi(x-x_0) Q(\sigma(x-x_0)),
\end{cases}
\end{equation}
where $A_\sigma > 0$ is chosen such that $(\phi_{1\sigma}, \phi_{2\sigma}) \in \mathcal{M}$, $\theta \in [0,1]$ and $x_0 \in \mathbb{R}^N$ will be determined later.
Since $(\phi_{1\sigma}, \phi_{2\sigma}) \in \mathcal{M}$, we have
\begin{equation}\label{2-19}
1 = \int_{\mathbb{R}^N} \frac{A_\sigma^2}{\|Q\|_2^2} \sigma^N \varphi^2(x-x_0) Q^2(\sigma(x-x_0)) dx = \int_{\mathbb{R}^N} \frac{A_\sigma^2}{\|Q\|_2^2} \varphi^2(\sigma^{-1}x) Q^2(x) dx.
\end{equation}
One can deduce from \eqref{decay} and \eqref{2-19} that
\begin{equation}
\begin{aligned}
	\left| 1 - A_\sigma^2 \right|
	&= \left| \int_{\mathbb{R}^N} \frac{A_\sigma^2}{\|Q \|_2^2} \left( \varphi^2\left( \sigma^{-1} x \right) - 1 \right)Q^2(x) dx \right| \\
	&\leq \int_{\mathbb{R}^N \setminus B_\sigma} \left| \frac{A_\sigma^2}{\| Q \|_2^2}  \left( \varphi^2\left( \sigma^{-1} x \right) - 1 \right)Q^2(x) \right| dx\\
	&\leq C \int_{\mathbb{R}^N \setminus B_\sigma} |x|^{-2N-4s} dx\\
	&\leq C \sigma^{-N-4s},
\end{aligned}
\end{equation}
for $\sigma$ sufficiently large, where \(B_\sigma := \{ x \in \mathbb{R}^N \mid |x| < \sigma \}\). Thus, we can deduce that
\begin{equation} \label{eq:A_sigma}
1 - O(\sigma^{-N-4s}) \leq A_\sigma^2 \leq 1 + O(\sigma^{-N-4s}) \quad \text{as } \sigma \to \infty.
\end{equation}
Now, let us turn to estimate each term in $E_{a_1,a_2,\beta}(\phi_{1\sigma}, \phi_{2\sigma})$. Denoting $\varphi_\sigma = \varphi(\sigma^{-1}x)$ and applying the scaling properties of the fractional Laplacian to the test function $\phi_{1\sigma}$ with scaling coefficient $\sigma$, we thus obtain
\begin{equation}
\int_{\mathbb{R}^N} |(-\Delta)^{\frac s2} \phi_{1\sigma}|^2 dx = \frac{\theta A_\sigma^2}{\|Q\|_2^2} \sigma^{2s} \int_{\mathbb{R}^N} |(-\Delta)^{\frac s2} (\varphi_\sigma Q)|^2 dx.
\end{equation}
By Lemma \ref{app:trunc-N1}, there exists \(C>0\) such that for all \(\sigma\) sufficiently large,
\begin{equation}\label{eq:trunc-est}
	\int_{\mathbb{R}^N}\big|(-\Delta)^{\frac s2}(\varphi_\sigma Q)\big|^2\,dx
	\le
	\int_{\mathbb{R}^N}\big|(-\Delta)^{\frac s2}Q\big|^2\,dx
	+C\sigma^{-4s}.
\end{equation}
One can infer from\eqref{eq:A_sigma}--\eqref{eq:trunc-est} and \eqref{eq:pohozaev_kinetic_mass} that
\begin{equation} \label{eq:kinetic1}
\begin{aligned}	
	\int_{\mathbb{R}^N} |(-\Delta)^{\frac s2} \phi_{1\sigma}|^2 dx &\le \frac{\theta A_\sigma^2}{\|Q\|_2^2} \sigma^{2s} \left( \int_{\mathbb{R}^N} |(-\Delta)^{\frac s2} Q|^2 dx + C\sigma^{-4s} \right) \\
	&= \frac{\theta A_\sigma^2}{\|Q\|_2^2} \sigma^{2s} \int_{\mathbb{R}^N} |(-\Delta)^{\frac s2} Q|^2 dx + \frac{C \theta A_\sigma^2}{\|Q\|_2^2} \sigma^{-2s}\\	
	&\le \theta\sigma^{2s}
	+O(\sigma^{-2s}).
\end{aligned}
\end{equation}
Apply the same steps to the test function $\phi_{2\sigma}$, one can also obtain
\begin{equation} \label{eq:kinetic2}
\int_{\mathbb{R}^N} |(-\Delta)^{\frac s2} \phi_{2\sigma}|^2 dx \leq (1-\theta) \sigma^{2s} + O(\sigma^{-2s}).
\end{equation}
For the nonlinear terms, we first claim
\begin{equation} \label{eq:nonlinear1}
\int_{\mathbb{R}^N} |\phi_{1\sigma}|^p dx \leq \frac{p \theta^{\frac p2} \sigma^{2s}}{2a^*} + O(\sigma^{-N-2s}).
\end{equation}
Indeed, by the scaling property, one can infer
\begin{equation}\label{phi1}
\int_{\mathbb{R}^N} |\phi_{1\sigma}|^p dx = \frac{\theta^{\frac p2} A_\sigma^p}{\|Q\|_2^p} \sigma^{N(\frac p2-1)} \int_{\mathbb{R}^N} \varphi^p(\sigma^{-1}x) |Q(x)|^p dx.
\end{equation}
From the decay of $Q$, we have
\begin{equation}\label{phi2}
\begin{aligned}
&\hphantom{=}\left| \int_{\mathbb{R}^N} \varphi^p(\sigma^{-1}x) |Q(x)|^p dx - \int_{\mathbb{R}^N} |Q(x)|^p dx \right| \\
&=\int_{\mathbb{R}^N} \big(1 - \varphi^p(\sigma^{-1} x)\big) |Q(x)|^p dx \\
&= \int_{|x|>\sigma} \big(1 - \varphi^p(\sigma^{-1} x)\big)|Q(x)|^p dx\\
&\leq \int_{|x|>\sigma} |Q(x)|^p dx\\
&\leq C\int_{|x|>\sigma} |x|^{p(-N-2s)} dx\\
&\leq C \sigma^{-p(N+2s)+N}\\
&\leq C \sigma^{-N-4s} ,
\end{aligned}
\end{equation}
for $\sigma$ sufficiently large. Combining \eqref{eq:A_sigma}, \eqref{phi2} and \eqref{pohozaev-eq}, one can deduce from \eqref{phi1} that
\begin{equation}
\begin{aligned}
\int_{\mathbb{R}^N} |\phi_{1\sigma}|^p dx
&\leq \frac{\theta^{\frac p2} \sigma^{2s}}{\|Q\|_2^p} \|Q\|_p^p + O(\sigma^{-N-2s})\\
&=\frac{p\theta^{\frac p2}\sigma^{2s}}{2a^*}  + O(\sigma^{-N-2s}).\\
\end{aligned}
\end{equation}
Thus, the claim \eqref{eq:nonlinear1} is proved. Note that, in a similar way, we can also obtain the following two estimates.
\begin{equation} \label{eq:nonlinear2}
\int_{\mathbb{R}^N} |\phi_{2\sigma}|^p dx \leq \frac{p (1-\theta)^{\frac p2} \sigma^{2s}}{2a^*} + O(\sigma^{-N-2s}),
\end{equation}
and
\begin{equation} \label{eq:nonlinear3}
\int_{\mathbb{R}^N} |\phi_{1\sigma}|^q |\phi_{2\sigma}|^q dx \leq \frac{p \theta^{\frac q2} (1-\theta)^{\frac q2} \sigma^{2s}}{2a^*} + O(\sigma^{-N-2s}).
\end{equation}
Furthermore, since the function $x \mapsto V_i(x) \varphi^2(x-x_0)$ ($i=1,2$) is bounded, by the dominated convergence theorem, we infer
\begin{equation} \label{eq:potential1}
\begin{aligned}
	\lim_{\sigma\to\infty}\int_{\mathbb{R}^N} V_1(x)\phi_{1\sigma}^2 dx
	&= \lim_{\sigma\to\infty}\int_{\mathbb{R}^N} \frac{\theta A_\sigma^2}{\|Q\|_2^2}\sigma^N V_1(x)\varphi^2(x-x_0) Q^2(\sigma(x-x_0)) dx \\
	&= \lim_{\sigma\to\infty}\int_{\mathbb{R}^N} \frac{\theta A_\sigma^2}{\|Q\|_2^2} V_1\left(\sigma^{-1}y + x_0\right)\varphi^2\left(\sigma^{-1}y\right) Q^2(y) dy \\
	&= \frac{\theta}{\|Q\|_2^2} V_1(x_0)\int_{\mathbb{R}^N} Q^2(y) dy = \theta V_1(x_0),
\end{aligned}
\end{equation}
and similarly we can also obtain
\begin{equation} \label{eq:potential2}
\lim_{\sigma \to \infty} \int_{\mathbb{R}^N} V_2(x) \phi_{2\sigma}^2 dx = (1-\theta) V_2(x_0).
\end{equation}
Then, it follows from \eqref{eq:energy}, \eqref{eq:kinetic1}-\eqref{eq:nonlinear1} and \eqref{eq:nonlinear2}-\eqref{eq:potential2} that
\begin{equation} \label{eq:test_energy_final}
\begin{aligned}
E_{a_1,a_2,\beta}(\phi_{1\sigma}, \phi_{2\sigma}) =& \sum_{i=1}^2 \int_{\mathbb{R}^N} \left( |(-\Delta)^{\frac s2} \phi_{i\sigma}|^2 + V_i(x) \phi_{i\sigma}^2 \right) dx \\
&- \frac{2}{p} \left[ a_1 \int_{\mathbb{R}^N} |\phi_{1\sigma}|^p dx + a_2 \int_{\mathbb{R}^N} |\phi_{2\sigma}|^p dx + 2\beta \int_{\mathbb{R}^N} |\phi_{1\sigma}|^q |\phi_{2\sigma}|^q dx \right] \\
\leq& \sigma^{2s} - \frac{a_1 \theta^q \sigma^{2s}}{a^*} - \frac{a_2 (1-\theta)^q \sigma^{2s}}{a^*} - \frac{2\beta \theta^{\frac q2} (1-\theta)^{\frac q2} \sigma^{2s}}{a^*} + O(\sigma^{-2s}) \\
&+ \int_{\mathbb{R}^N} V_1(x) \phi_{1\sigma}^2 dx + \int_{\mathbb{R}^N} V_2(x) \phi_{2\sigma}^2 dx \\
=& \frac{\sigma^{2s}}{a^*} \left[ a^* - a_1 \theta^q - a_2 (1-\theta)^q - 2\beta \theta^{\frac q2} (1-\theta)^{\frac q2} \right] + O(\sigma^{-2s}) \\
&+ \theta V_1(x_0) + (1-\theta) V_2(x_0).
\end{aligned}
\end{equation}
Next, we finish the proof of (ii) by considering the following three cases.

\noindent\textbf{Case 1: $a_1 > a^*$.} One can obtain from \eqref{eq:test_energy_final} by taking $\theta = 1$ and any $x_0 \in \mathbb{R}^N$ that
\begin{equation}
E_{a_1,a_2,\beta}(\phi_{1\sigma}, \phi_{2\sigma}) \leq \frac{\sigma^{2s}}{a^*} (a^* - a_1) + V_1(x_0) + O(\sigma^{-2s}) \to -\infty \quad \text{as } \sigma \to \infty,
\end{equation}
which implies $e(a_1,a_2,\beta) = -\infty$ in this case. Therefore, if $a_1 > a^*$, we deduce that there is no minimizer of \eqref{eq:min_problem}.

\noindent\textbf{Case 2: $a_2 > a^*$.} By taking $\theta = 0$ and any $x_0 \in \mathbb{R}^N$ in \eqref{eq:test_energy_final}, one can infer
\begin{equation}
E_{a_1,a_2,\beta}(\phi_{1\sigma}, \phi_{2\sigma}) \leq \frac{\sigma^{2s}}{a^*} (a^* - a_2) + V_2(x_0) + O(\sigma^{-2s}) \to -\infty \quad \text{as } \sigma \to \infty,
\end{equation}
which indicates $e(a_1,a_2,\beta) = -\infty$. Hence, if $a_2 > a^*$, we conclude that there is no minimizer of \eqref{eq:min_problem}.

\noindent\textbf{Case 3: $\beta > \beta^*$.} Without loss of generality, one can suppose $0 \leq a_1, a_2 < a^*$ in this case. From the definition of $\beta^*$, there exists some $\theta_0 \in (0,1)$ such that
\begin{equation}
a^* - a_1 \theta_0^q - a_2 (1-\theta_0)^q = 2\beta^* \theta_0^{\frac q2} (1-\theta_0)^{\frac q2}.
\end{equation}

Taking $\theta = \theta_0$ and any $x_0 \in \mathbb{R}^N$ in \eqref{eq:test_energy_final}, we deduce that
\begin{equation}
\begin{aligned}
E_{a_1,a_2,\beta}(\phi_{1\sigma}, \phi_{2\sigma}) \leq& \frac{2\sigma^{2s}}{a^*} \theta_0^{\frac q2} (1-\theta_0)^{\frac q2} (\beta^* - \beta) + \theta_0 V_1(x_0) + (1-\theta_0) V_2(x_0) \\
&+ O(\sigma^{-2s}) \to -\infty \quad \text{as } \sigma \to \infty,
\end{aligned}
\end{equation}
which shows $e(a_1,a_2,\beta) = -\infty$ as well. We thus conclude that, if $\beta > \beta^*$, there is no minimizer of \eqref{eq:min_problem}.

(iii) Finally, we prove that there is no minimizer of \eqref{eq:min_problem} in the case $a_1 = a^*$ or $a_2 = a^*$ or $\beta = \beta^*$. By (ii), it remains to consider the parameters satisfying $0 \le a_1, a_2 \le a^*$ and $\Lambda(a_1, a_2, \beta) \le a^*$; for these parameters, \eqref{eq:lower_bound} yields $E_{a_1, a_2, \beta}(u_1, u_2) \ge 0$ for every $(u_1, u_2) \in \mathcal{M}$.
From the assumption $(\Omega)$ we know that
\begin{equation}
\mathcal{V} = \left\{ x \in \mathbb{R}^N : V_1(x) = V_2(x) = 0 \right\} \neq \emptyset.
\end{equation}
Therefore, we can choose $x_0 \in \mathcal{V}$ such that $V_1(x_0) = V_2(x_0) = 0$. It then follows from \eqref{eq:lower_bound} and \eqref{eq:test_energy_final} that
\begin{equation} \label{eq:critical_energy}
\begin{aligned}
0 \leq e(a_1,a_2,\beta) \leq {}& E_{a_1,a_2,\beta}(\phi_{1\sigma}, \phi_{2\sigma})\\
\leq {}& \frac{\sigma^{2s}}{a^*} \Bigl[ a^* - a_1 \theta^q - a_2 (1-\theta)^q - 2\beta \theta^{q/2} (1-\theta)^{q/2} \Bigr] + O(\sigma^{-2s}).
\end{aligned}
\end{equation}
We next discuss the following three cases similarly as part (ii).

\noindent\textbf{Case 1:} $a_1 = a^*$. Let $\theta = 1$ in \eqref{eq:critical_energy} and $\sigma \to \infty$, we then obtain
\begin{equation}\label{hhh}
0 \leq e(a^*, a_2, \beta) \leq \lim_{\sigma \to \infty} E_{a^*, a_2, \beta}(\phi_{1\sigma}, \phi_{2\sigma}) = 0.
\end{equation}
Next, we argue by contradiction. Suppose $e(a^*, a_2, \beta)$ has a minimizer $(u_{10}, u_{20})$. From \eqref{eq:nonlinear}, \eqref{eq:lower_bound} and \eqref{hhh}, we infer that
\begin{equation}\label{eq}
\int_{\mathbb{R}^N} \left( |(-\Delta)^{s/2} u_{10}|^2 + |(-\Delta)^{s/2} u_{20}|^2 \right) dx = \frac{2a^*}{p} \int_{\mathbb{R}^N} \bigl(u_{10}^2 + u_{20}^2\bigr)^{q} dx,
\end{equation}
\begin{equation}\label{c}
\int_{\mathbb{R}^N} V_1(x) u_{10}^2 dx = \int_{\mathbb{R}^N} V_2(x) u_{20}^2 dx = 0.
\end{equation}
One can deduce from Lemma \ref{lem:optimizers} and the equality \eqref{eq} that
\begin{equation}\label{cc}
(u_{10}(x), u_{20}(x)) = (\tau \sin\theta \, Q(\alpha x + \eta), \tau \cos\theta \, Q(\alpha x + \eta))
\end{equation}
with some $\tau>0, \alpha>0, \theta\in[0,2\pi)$ and $\eta\in\mathbb{R}^N$.
 However, \eqref{c} shows that the support of $(u_{10},u_{20})$ is compact, which contradicts \eqref{cc}. Therefore, in this case, there is no minimizer of $e(a^*,a_2,\beta)$.

\noindent\textbf{Case 2:} $a_2 = a^*$. Let $\theta = 0$ in \eqref{eq:critical_energy} and $\sigma \to \infty$, we then have
\begin{equation}
0 \leq e(a_1, a^*, \beta) \leq \lim_{\sigma \to \infty} E_{a_1, a^*, \beta}(\phi_{1\sigma}, \phi_{2\sigma}) = 0,
\end{equation}
which shows that there is no minimizer of $e(a_1, a^*, \beta)$ by arguing as above.

\noindent\textbf{Case 3:} $\beta = \beta^*$ and $0 \leq a_1, a_2  < a^*$. Let $\theta = \theta_0$ in \eqref{eq:critical_energy} and $\sigma \to \infty$, we then get
\begin{equation}\label{ss}
0 \leq e(a_1, a_2, \beta^*) \leq \lim_{\sigma \to \infty} E_{a_1, a_2, \beta^*}(\phi_{1\sigma}, \phi_{2\sigma}) = 0,
\end{equation}
which implies that there is no minimizer of $e(a_1, a_2, \beta^*)$ by arguing as above.

Furthermore, in the case $0 \leq a_1, a_2 < a^*$ and $0 < \beta < \beta^*$, let $\theta = \theta_0$ and $\sigma = (\beta^* - \beta)^{-\frac{1}{4s}}$ in \eqref{eq:critical_energy}, one can infer
\begin{equation}
0 \leq e(a_1,a_2,\beta) \leq \frac{\sigma^{2s}}{a^*} \cdot 2 \theta_0^{\frac q2} (1-\theta_0)^{\frac q2} (\beta^* - \beta) + O(\sigma^{-2s}) \to 0 \quad \text{as } \beta \to \beta^*.
\end{equation}
Combining this and \eqref{ss}, we thus deduce that
\begin{equation}
\lim_{\beta \to \beta^*} e(a_1,a_2,\beta) = e(a_1,a_2,\beta^*) = 0, \quad 0 \leq a_1, a_2 < a^*,
\end{equation}
which completes the proof of Theorem \ref{thm:main}.
\end{proofof}

\section{General concentration behavior}
In this section, we investigate the concentration behavior of the minimizers for the problem \eqref{eq:min_problem} under general trapping potentials as $\beta \to \beta^*$. We specifically focus on the one-dimensional, half-order fractional Laplacian regime, which implies $p=4$ and $q=2$. Since \eqref{convex} gives $\|(-\Delta)^{s/2}|u_i|\|_2\le\|(-\Delta)^{s/2}u_i\|_2$ and $E_{a_1,a_2,\beta}$ depends on $u_i$ only through $|u_i|$, we may assume that the minimizers are nonnegative. For this purpose, we start with the following two lemmas.

\begin{lemma}
	Suppose $(\Omega)$ holds and let $(u_{1\beta},u_{2\beta})$ be a nonnegative minimizer of $e(a_1,a_2,\beta)$ for $0\le a_1,a_2<a^*$ and $0<\beta<\beta^*$. Then the following hold:
	\begin{enumerate}[label=(\roman*)]
		\item $(u_{1\beta},u_{2\beta})$ blows up as $\beta\to\beta^*$, i.e., for $i=1,2$,
		\begin{equation}\label{eq:3-2}
			\lim_{\beta \to \beta^*}\int_{\mathbb R} |(-\Delta)^{\frac14}u_{i\beta}|^2 dx=\infty \quad\text{and}\quad \lim_{\beta\to\beta^*}\int_{\mathbb R} u_{i\beta}^4 dx=\infty.
		\end{equation}
		\item $(u_{1\beta},u_{2\beta})$ also satisfies
		\begin{equation}\label{eq:3-3}
			\lim_{\beta\to\beta^*}\int_{\mathbb R} V_1(x)u_{1\beta}^2 dx =\lim_{\beta\to\beta^*}\int_{\mathbb R} V_2(x)u_{2\beta}^2 dx=0,
		\end{equation}
		\begin{equation}\label{eq:3-4}
			\lim_{\beta\to\beta^*}\Big[(\beta^*-\beta)\int_{\mathbb R} u_{1\beta}^2 u_{2\beta}^2 dx\Big]=0,
		\end{equation}
		\begin{equation}\label{eq:3-5}
			\lim_{\beta\to\beta^*}\Big[\sqrt{a^*-a_1}\sqrt{\int_{\mathbb R}u_{1\beta}^4 dx} -\sqrt{a^*-a_2}\sqrt{\int_{\mathbb R}u_{2\beta}^4 dx}\Big]^2=0,
		\end{equation}
		\begin{equation}\label{eq:3-6}
			\lim_{\beta\to\beta^*}\Big[\sqrt{\int_{\mathbb R}u_{1\beta}^4 dx}\sqrt{\int_{\mathbb R}u_{2\beta}^4 dx} -\int_{\mathbb R}u_{1\beta}^2 u_{2\beta}^2 dx\Big]=0,
		\end{equation}
		as well as
		\begin{equation}\label{eq:3-7}
			\lim_{\beta\to\beta^*} \frac{\int_{\mathbb R}\big(|(-\Delta)^{\frac14}u_{1\beta}|^2 + |(-\Delta)^{\frac14}u_{2\beta}|^2\big)dx} {\int_{\mathbb R}(u_{1\beta}^2+u_{2\beta}^2)^2 dx} =\frac{a^*}{2},
		\end{equation}
		\begin{equation}\label{eq:3-8}
			\lim_{\beta\to\beta^*}\frac{\int_{\mathbb R}u_{1\beta}^4 dx}{\int_{\mathbb R}u_{2\beta}^4 dx} =\frac{a^*-a_2}{a^*-a_1}.
		\end{equation}
	\end{enumerate}
\end{lemma}

\begin{proof}
	Combining the refined Gagliardo--Nirenberg inequality \eqref{refine} and the Schwarz inequality, we deduce that
\begin{equation}\label{eq:3-9}
\begin{aligned}
	E_{a_1,a_2,\beta}(u_{1\beta},u_{2\beta})
	&=\sum_{i=1}^2\int_{\mathbb R}|(-\Delta)^{\frac14}u_{i\beta}|^2 dx -\frac{a^*}{2}\int_{\mathbb R}(u_{1\beta}^2+u_{2\beta}^2)^2 dx \\
	&\quad+\frac12\Big[\sqrt{a^*-a_1}\sqrt{\int_{\mathbb R}u_{1\beta}^4 dx} -\sqrt{a^*-a_2}\sqrt{\int_{\mathbb R}u_{2\beta}^4 dx}\Big]^2 \\
	&\quad+\sqrt{(a^*-a_1)(a^*-a_2)} \Big[\sqrt{\int_{\mathbb R}u_{1\beta}^4 dx}\sqrt{\int_{\mathbb R}u_{2\beta}^4 dx} -\int_{\mathbb R}u_{1\beta}^2u_{2\beta}^2 dx\Big] \\
	&\quad+(\beta^*-\beta)\int_{\mathbb R}u_{1\beta}^2u_{2\beta}^2 dx +\sum_{i=1}^2\int_{\mathbb R}V_i(x)u_{i\beta}^2 dx\\
	&\ge 0.
\end{aligned}
\end{equation}
One can immediately deduce from \eqref{eq:energy_limit} that \eqref{eq:3-3}--\eqref{eq:3-6} and the following hold,
\begin{equation}\label{eq:3-10}
	\lim_{\beta\to\beta^*}\Big[\int_{\mathbb R}\big(|(-\Delta)^{\frac14}u_{1\beta}|^2 + |(-\Delta)^{\frac14}u_{2\beta}|^2\big) -\frac{a^*}{2}\int_{\mathbb R}(u_{1\beta}^2+u_{2\beta}^2)^2 dx\Big]=0.
\end{equation}
Note that
\begin{equation}\label{eq:3-11}
	\lim_{\beta\to\beta^*}\int_{\mathbb R}\big(|(-\Delta)^{\frac14}u_{1\beta}|^2 + |(-\Delta)^{\frac14}u_{2\beta}|^2\big)dx=\infty.
\end{equation}
In fact, suppose this does not hold, then we can infer that there exists a sequence $\{\beta_k\}$ satisfying $\beta_k\to\beta^*$ as $k\to\infty$, such that $\{(u_{1\beta_k},u_{2\beta_k})\}$ is bounded uniformly in $\mathcal{H}_1\times\mathcal{H}_2$. Therefore, from Lemma \ref{lem:compact_embedding}, we deduce that there exists a subsequence of $\{\beta_k\}$, still denoted by $\{\beta_k\}$, and $(u_{10},u_{20})\in \mathcal{H}_1\times\mathcal{H}_2$ such that
\begin{equation}\label{eq:3-12}
	\begin{cases}
		(u_{1\beta_k},u_{2\beta_k})\rightharpoonup (u_{10},u_{20}) \quad\text{weakly in } \mathcal{H}_1\times\mathcal{H}_2,\\
		(u_{1\beta_k},u_{2\beta_k})\to (u_{10},u_{20}) \quad\text{strongly in } L^r(\mathbb R)\times L^r(\mathbb R) \text{ for } 2\le r<\infty.
	\end{cases}
\end{equation}
Then, it follows from Lemma \ref{lem:convergence} and the weak lower semi-continuity of the norm $\|\cdot\|_{\mathcal{H}_i}$ that
\begin{equation}\label{eq:3-13}
	0=e(a_1,a_2,\beta^*)\le E_{a_1,a_2,\beta^*}(u_{10},u_{20}) \le \lim_{k\to\infty}E_{a_1,a_2,\beta_k}(u_{1\beta_k},u_{2\beta_k}) =\lim_{k\to\infty}e(a_1,a_2,\beta_k)=0.
\end{equation}
It follows that $(u_{10},u_{20})$ is a minimizer of $E_{a_1,a_2,\beta^*}$, which contradicts Theorem \ref{thm:main} and hence \eqref{eq:3-11} holds. Combining \eqref{eq:3-10} and \eqref{eq:3-11}, \eqref{eq:3-7} thus follows. Moreover, from \eqref{eq:3-6}, \eqref{eq:3-7} and \eqref{eq:3-11}, we further conclude that
\begin{equation}\label{eq:3-14}
	\lim_{\beta\to\beta^*}\int_{\mathbb R}(u_{1\beta}^2+u_{2\beta}^2)^2 dx =\lim_{\beta\to\beta^*}\Big(\sqrt{\int_{\mathbb R}u_{1\beta}^4 dx}+\sqrt{\int_{\mathbb R}u_{2\beta}^4 dx}\Big)^2=\infty.
\end{equation}
It then follows from \eqref{eq:3-5} that for $i=1,2$,
\begin{equation}\label{eq:3-15}
	\lim_{\beta\nearrow\beta^*}\int_{\mathbb R}u_{i\beta}^4 dx=\infty.
\end{equation}
Therefore, we can deduce from the Gagliardo--Nirenberg inequality \eqref{lem:single_gns} that \eqref{eq:3-2} holds. Finally, \eqref{eq:3-2} and \eqref{eq:3-5} show that \eqref{eq:3-8} is true, which completes the proof.
\end{proof}

For the nonnegative minimizer $(u_{1\beta},u_{2\beta})$ of \eqref{eq:min_problem}, define
\begin{equation}\label{eq:3-16}
	\varepsilon_{\beta}:=\Big(\int_{\mathbb R}\big(|(-\Delta)^{\frac14}u_{1\beta}|^2 + |(-\Delta)^{\frac14}u_{2\beta}|^2\big)dx\Big)^{-1}>0,
\end{equation}
which satisfies $\varepsilon_{\beta}\to0$ as $\beta\to\beta^*$ according to \eqref{eq:3-11}. Next, we introduce another lemma that will play a fundamental role in the proof of Theorem \ref{con}.

\begin{lemma}\label{lemma3-2}
	Suppose $(\Omega)$ holds and let $(u_{1\beta},u_{2\beta})$ be a nonnegative minimizer of $e(a_1,a_2,\beta)$ for $0\le a_1,a_2<a^*$ and $0<\beta<\beta^*$. Then for any sequence $\beta_k\to\beta^*$ as $k\to\infty$, we have
	\begin{enumerate}[label=(\roman*)]
		\item There exist a sequence $\{\eta_{\beta_k}\}\subset\mathbb R$ and constants $r_0,\delta>0$ such that
		\begin{equation}\label{eq:3-17}
			w_{i\beta_k}(x):=\varepsilon_{\beta_k}^{\frac12}u_{i\beta_k}\big(\varepsilon_{\beta_k}x+\varepsilon_{\beta_k}\eta_{\beta_k}\big)
		\end{equation}
		satisfies
		\begin{equation}\label{eq:3-18}
			\liminf_{k\to\infty}\int_{B_{r_0}(0)}|w_{i\beta_k}(x)|^2 dx\ge \delta>0,\quad i=1,2.
		\end{equation}
		Furthermore,
		\begin{equation}\label{eq:3-19}
			\int_{\mathbb R}\big(|(-\Delta)^{\frac14}w_{1\beta_k}|^2 + |(-\Delta)^{\frac14}w_{2\beta_k}|^2\big)dx =\int_{\mathbb R}(w_{1\beta_k}^2+w_{2\beta_k}^2)dx=1,
		\end{equation}
		\begin{equation}\label{eq:3-20}
			\lim_{k\to\infty}\biggl[\frac{\beta^*-\beta_k}{\varepsilon_{\beta_k}}\int_{\mathbb R}w_{1\beta_k}^2 w_{2\beta_k}^2 dx\biggr]=0,
		\end{equation}
		\begin{equation}\label{eq:3-21}
			\lim_{k\to\infty}\int_{\mathbb R} \big(w_{1\beta_k}^2 + w_{2\beta_k}^2\big)^2 dx = \frac{2}{a^*},
		\end{equation}
		\begin{equation}\label{eq:3-22}
			\begin{aligned}
&\lim_{k\to\infty}\int_{\mathbb R} w_{1\beta_k}^4 dx = \frac{2\gamma^2}{a^*},\quad
\lim_{k\to\infty}\int_{\mathbb R} w_{2\beta_k}^4 dx = \frac{2(1-\gamma)^2}{a^*},\\
&\lim_{k\to\infty}\int_{\mathbb R} w_{1\beta_k}^2 w_{2\beta_k}^2 dx = \frac{2\gamma(1-\gamma)}{a^*}.
\end{aligned}
		\end{equation}
		\item There exists a subsequence of $\{\beta_k\}$, still denoted by $\{\beta_k\}$, such that
		\begin{equation}\label{eq:3-23}
			x_{\beta_k}:=\varepsilon_{\beta_k}\eta_{\beta_k} \to x_0 \quad\text{as }k\to\infty
		\end{equation}
		for some $x_0\in\mathcal{V}$.
		\item There exists a subsequence of $\{\beta_k\}$, still denoted by $\{\beta_k\}$, and $\eta_0\in\mathbb R$ such that
		\begin{equation}\label{eq:3-24}
			\begin{cases}
				\displaystyle\lim_{k\to\infty} w_{1\beta_k} = \frac{\sqrt{\gamma}}{\|Q\|_2} Q(x-\eta_0),\\
				\displaystyle\lim_{k\to\infty} w_{2\beta_k} = \frac{\sqrt{1-\gamma}}{\|Q\|_2} Q(x-\eta_0),
			\end{cases}
		\end{equation}
		strongly in $H^{\frac12}(\mathbb R)$.
	\end{enumerate}
\end{lemma}

\begin{proof}
(i) Denote
\begin{equation}\label{eq:3-25}
	v_{i\beta_k}:=\varepsilon_{\beta_k}^{\frac12}u_{i\beta_k}\big(\varepsilon_{\beta_k}x\big),\quad i=1,2.
\end{equation}
By a simple calculation, one obtains from \eqref{eq:3-16} and the mass constraint that
\begin{equation}\label{eq:3-26}
	\begin{aligned}
		&\hphantom{==}\int_{\mathbb R}\big(|(-\Delta)^{\frac14}v_{1\beta_k}|^2 + |(-\Delta)^{\frac14}v_{2\beta_k}|^2\big)dx\\
		&=\varepsilon_{\beta_k}\int_{\mathbb R}\big(|(-\Delta)^{\frac14}u_{1\beta_k}|^2 + |(-\Delta)^{\frac14}u_{2\beta_k}|^2\big)dx\\
		&=1\\
		&=\int_{\mathbb R}\big(|u_{1\beta_k}|^2 + |u_{2\beta_k}|^2\big)dx\\
		&=\int_{\mathbb R}\big(|v_{1\beta_k}|^2 + |v_{2\beta_k}|^2\big)dx.
	\end{aligned}
\end{equation}
Besides, it comes from $\int_{\mathbb R}v_{1\beta_k}^2v_{2\beta_k}^2dx=\varepsilon_{\beta_k}\int_{\mathbb R}u_{1\beta_k}^2u_{2\beta_k}^2dx$ and \eqref{eq:3-4} that
\begin{equation}\label{eq:3-27}
	\lim_{k\to\infty}\Big[\frac{\beta^*-\beta_k}{\varepsilon_{\beta_k}}\int_{\mathbb R}v_{1\beta_k}^2v_{2\beta_k}^2dx\Big]
	=\lim_{k\to\infty}\big[(\beta^*-\beta_k)\int_{\mathbb R}u_{1\beta_k}^2u_{2\beta_k}^2dx\big]=0.
\end{equation}
Similarly from \eqref{eq:3-5} and \eqref{eq:3-6}, we can obtain
\begin{equation}\label{eq:3-28}
	\lim_{k\to\infty}\Big[\sqrt{a^*-a_1}\sqrt{\int_{\mathbb R}v_{1\beta_k}^4dx} -\sqrt{a^*-a_2}\sqrt{\int_{\mathbb R}v_{2\beta_k}^4dx}\Big]^2=0
\end{equation}
and
\begin{equation}\label{eq:3-29}
	\lim_{k\to\infty}\Big[\sqrt{\int_{\mathbb R}v_{1\beta_k}^4dx}\sqrt{\int_{\mathbb R}v_{2\beta_k}^4dx} -\int_{\mathbb R}v_{1\beta_k}^2v_{2\beta_k}^2dx\Big]=0.
\end{equation}
Moreover, from \eqref{eq:3-7} and \eqref{eq:3-26}, one can easily conclude that
\begin{equation}\label{eq:3-30}
	\begin{aligned}
		\lim_{k\to\infty}\int_{\mathbb R}\big(v_{1\beta_k}^2+v_{2\beta_k}^2\big)^2dx
		&=\lim_{k\to\infty}\frac{\int_{\mathbb R}\big(v_{1\beta_k}^2+v_{2\beta_k}^2\big)^2dx} {\int_{\mathbb R}\big(|(-\Delta)^{\frac14}v_{1\beta_k}|^2+|(-\Delta)^{\frac14}v_{2\beta_k}|^2\big)dx}\\
		&=\lim_{k\to\infty}\frac{\int_{\mathbb R}\big(u_{1\beta_k}^2+u_{2\beta_k}^2\big)^2dx} {\int_{\mathbb R}\big(|(-\Delta)^{\frac14}u_{1\beta_k}|^2+|(-\Delta)^{\frac14}u_{2\beta_k}|^2\big)dx}\\
		&=\frac{2}{a^*}.
	\end{aligned}
\end{equation}
Combining \eqref{eq:3-28}--\eqref{eq:3-30} and \eqref{eq:3-1}, we can deduce that
\begin{equation}\label{eq:3-31}
	\begin{aligned}
		\frac{2}{a^*} &= \lim_{k\to\infty}\int_{\mathbb R}\big(v_{1\beta_k}^2+v_{2\beta_k}^2\big)^2dx\\
		&=\lim_{k\to\infty}\Big(\sqrt{\int_{\mathbb R}v_{1\beta_k}^4dx}+\sqrt{\int_{\mathbb R}v_{2\beta_k}^4dx}\Big)^2\\
		&=\frac{1}{\gamma^2}\lim_{k\to\infty}\int_{\mathbb R}v_{1\beta_k}^4dx\\
		&=\frac{1}{(1-\gamma)^2}\lim_{k\to\infty}\int_{\mathbb R}v_{2\beta_k}^4dx,
	\end{aligned}
\end{equation}
where $\gamma$ is defined by \eqref{eq:3-1}. Therefore, we have
\begin{equation}\label{eq:3-32}
	\lim_{k\to\infty}\int_{\mathbb R}v_{1\beta_k}^4dx=\frac{2\gamma^2}{a^*},\quad \lim_{k\to\infty}\int_{\mathbb R}v_{2\beta_k}^4dx=\frac{2(1-\gamma)^2}{a^*},\quad \lim_{k\to\infty}\int_{\mathbb R}v_{1\beta_k}^2v_{2\beta_k}^2dx=\frac{2\gamma(1-\gamma)}{a^*}.
\end{equation}
Note from \eqref{eq:3-26} that $\{v_{i\beta_k}\}$ is bounded in $H^{\frac12}(\mathbb R)$ for $i=1,2$. Next, we claim that there are a sequence $\{\eta_{\beta_k}\}\subset\mathbb{R}$ and constants $r_0,\delta>0$ such that
\begin{equation}\label{eq:3-33}
	\liminf_{k\to\infty}\int_{B_{r_0}(\eta_{\beta_k})}|v_{1\beta_k}(x)|^2dx\ge \delta>0.
\end{equation}
Arguing by contradiction, suppose that \eqref{eq:3-33} is false. Then, for any $r>0$, there exists a subsequence of $\{v_{1\beta_k}\}$, still denoted by $\{v_{1\beta_k}\}$, such that
\begin{equation}\label{eq:3-34}
	\limsup_{k\to\infty}\sup_{y\in\mathbb R}\int_{B_r(y)}|v_{1\beta_k}(x)|^2dx=0,
\end{equation}
from which we can deduce that $v_{1\beta_k}\to0$ in $L^m(\mathbb R)$ for all $2<m<\infty$ by Lions' vanishing lemma (see \cite[Lemma 2.4]{GuoWangZengZhou2018}). However, this implies $\int_{\mathbb R}v_{1\beta_k}^4dx\to0$ as $k\to\infty$, which sharply contradicts \eqref{eq:3-32}. Therefore, \eqref{eq:3-33} follows. From \eqref{eq:3-17}, we thus deduce
\begin{equation}\label{eq:3-35}
	\liminf_{k\to\infty}\int_{B_{r_0}(0)}|w_{1\beta_k}(x)|^2dx\ge \delta>0.
\end{equation}
Comparing \eqref{eq:3-17} with \eqref{eq:3-25}, it follows from \eqref{eq:3-26}--\eqref{eq:3-32} that \eqref{eq:3-19}--\eqref{eq:3-22} hold. Furthermore, from \eqref{eq:3-19}, one can assume that $w_{i\beta_k}\rightharpoonup w_{i0}\ge0$ in $H^{\frac12}(\mathbb R)$ for some $w_{i0}\in H^{\frac12}(\mathbb R)$, where $i=1,2$. Moreover, from \eqref{eq:3-35}, one can easily see that $w_{10}\neq0$.

 To complete the proof, it remains to obtain a result similar to \eqref{eq:3-35} for $w_{2\beta_k}$ with the same sequence $\eta_{\beta_k}$. We first claim that $\{x_{\beta_k}\}$ defined in \eqref{eq:3-23} is bounded in $\mathbb R$. In fact, if it does not hold, there exists a subsequence of $\{x_{\beta_k}\}$, still denoted by $\{x_{\beta_k}\}$, such that $|x_{\beta_k}|\to\infty$ as $k\to\infty$. From \eqref{eq:3-3}, we infer
\begin{equation}\label{eq:3-36}
	\sum_{i=1}^2\int_{\mathbb R}V_i\big(\varepsilon_{\beta_k}x+x_{\beta_k}\big)|w_{i\beta_k}(x)|^2dx
	=\sum_{i=1}^2\int_{\mathbb R}V_i(x)|u_{i\beta_k}(x)|^2dx\to0 \quad\text{as }k\to\infty.
\end{equation}
Since $V_1(x)\to\infty$ as $|x|\to\infty$, for any $C>0$ and sufficiently large $k$, we have $V_1(\varepsilon_{\beta_k}x+x_{\beta_k})\ge C$ on $B_{r_0}(0)$. Combined with \eqref{eq:3-35}, this yields
\[
\liminf_{k\to\infty}\sum_{i=1}^2\int_{\mathbb R}V_i\big(\varepsilon_{\beta_k}x+x_{\beta_k}\big)|w_{i\beta_k}(x)|^2dx
\ge C\liminf_{k\to\infty}\int_{B_{r_0}(0)}|w_{1\beta_k}(x)|^2dx\ge C \delta>0,
\]
which contradicts \eqref{eq:3-36} and therefore the claim follows. Then, we can infer that, going if necessary to a subsequence, there exists some point $x_0\in\mathbb R$ such that
\begin{equation}\label{eq:3-37}
	x_{\beta_k}\to x_0 \quad\text{as }k\to\infty.
\end{equation}
Applying the translation and scaling properties of the fractional Laplacian $(-\Delta)^{\frac12}$ to $w_{i\beta_k}$, one can deduce from \eqref{eq:3-17} that for $i=1,2$,
\begin{equation}\label{eq:3-38}
	(-\Delta)^{\frac12}w_{i\beta_k}(x) = \varepsilon_{\beta_k}^{\frac32}(-\Delta)^{\frac12}u_{i\beta_k}\big(\varepsilon_{\beta_k}x+x_{\beta_k}\big).
\end{equation}
Moreover, from \eqref{eq:3-17}, one can easily obtain
\begin{equation}\label{eq:3-39}
	w_{1\beta_k}^3(x) = \varepsilon_{\beta_k}^{\frac32}u_{1\beta_k}^3\big(\varepsilon_{\beta_k}x+x_{\beta_k}\big),
\end{equation}
\begin{equation}\label{eq:3-40}
	w_{2\beta_k}^2(x)w_{1\beta_k}(x) = \varepsilon_{\beta_k}^{\frac32}u_{2\beta_k}^2\big(\varepsilon_{\beta_k}x+x_{\beta_k}\big)u_{1\beta_k}\big(\varepsilon_{\beta_k}x+x_{\beta_k}\big),
\end{equation}
\begin{equation}\label{eq:3-41}
	w_{1\beta_k}^2(x)w_{2\beta_k}(x) = \varepsilon_{\beta_k}^{\frac32}u_{1\beta_k}^2\big(\varepsilon_{\beta_k}x+x_{\beta_k}\big)u_{2\beta_k}\big(\varepsilon_{\beta_k}x+x_{\beta_k}\big).
\end{equation}
Since the nonnegative minimizer $(u_{1\beta_k},u_{2\beta_k})$ satisfies system \eqref{eq:sys} with the Lagrange multiplier $\mu_{\beta_k}\in\mathbb R$ for $p=4$ and $q=2$, we then conclude from \eqref{eq:3-17} and \eqref{eq:3-38}--\eqref{eq:3-41} that $(w_{1\beta_k},w_{2\beta_k})$ satisfies the following equations in $\mathbb R$,
\begin{equation}\label{eq:3-42}
	\begin{cases}
		(-\Delta)^{\frac12}w_{1\beta_k} + \varepsilon_{\beta_k}V_1\big(\varepsilon_{\beta_k}x+x_{\beta_k}\big)w_{1\beta_k} = \varepsilon_{\beta_k}\mu_{\beta_k}w_{1\beta_k} + a_1w_{1\beta_k}^3 + \beta_k w_{2\beta_k}^2 w_{1\beta_k},\\
		(-\Delta)^{\frac12}w_{2\beta_k} + \varepsilon_{\beta_k}V_2\big(\varepsilon_{\beta_k}x+x_{\beta_k}\big)w_{2\beta_k} = \varepsilon_{\beta_k}\mu_{\beta_k}w_{2\beta_k} + a_2w_{2\beta_k}^3 + \beta_k w_{1\beta_k}^2 w_{2\beta_k}.
	\end{cases}
\end{equation}
From \eqref{eq:3-42} and \eqref{eq:3-19}, one obtains that
\begin{equation}\label{eq:3-43}
	\begin{aligned}
		\varepsilon_{\beta_k}\mu_{\beta_k} &= \sum_{i=1}^2\int_{\mathbb R}\big(|(-\Delta)^{\frac14}w_{i\beta_k}|^2 + \varepsilon_{\beta_k}V_i(\varepsilon_{\beta_k}x+x_{\beta_k})w_{i\beta_k}^2\big)dx\\
		&\quad - \big(a_1\int_{\mathbb R}w_{1\beta_k}^4dx + a_2\int_{\mathbb R}w_{2\beta_k}^4dx + 2\beta_k\int_{\mathbb R}w_{1\beta_k}^2w_{2\beta_k}^2dx\big).
	\end{aligned}
\end{equation}
On the other hand, one can deduce from \eqref{eq:energy} and \eqref{eq:3-38}--\eqref{eq:3-41} that
\begin{equation}\label{eq:3-44}
	\begin{aligned}
		\varepsilon_{\beta_k}e(a_1,a_2,\beta_k)
		&=\varepsilon_{\beta_k}E_{a_1,a_2,\beta_k}(u_{1\beta_k},u_{2\beta_k})\\
		&=\sum_{i=1}^2\int_{\mathbb R}\big(|(-\Delta)^{\frac14}w_{i\beta_k}|^2+\varepsilon_{\beta_k}V_i(\varepsilon_{\beta_k}x+x_{\beta_k})w_{i\beta_k}^2\big)dx\\
		&\quad -\frac12\big(a_1\int_{\mathbb R}w_{1\beta_k}^4dx+a_2\int_{\mathbb R}w_{2\beta_k}^4dx+2\beta_k\int_{\mathbb R}w_{1\beta_k}^2w_{2\beta_k}^2dx\big).
	\end{aligned}
\end{equation}
Combining \eqref{eq:3-43} and \eqref{eq:3-44}, we conclude from \eqref{eq:3-16} and \eqref{eq:3-19} that
\begin{equation}\label{eq:3-45}
	\begin{aligned}
		\varepsilon_{\beta_k}\mu_{\beta_k} &= 2\varepsilon_{\beta_k}e(a_1,a_2,\beta_k) - \sum_{i=1}^2\int_{\mathbb R}\big(|(-\Delta)^{\frac14}w_{i\beta_k}|^2+\varepsilon_{\beta_k}V_i(\varepsilon_{\beta_k}x+x_{\beta_k})w_{i\beta_k}^2\big)dx\\
		&\to-1 \quad \text{as } k\to\infty.
	\end{aligned}
\end{equation}
Since the right-hand side of \eqref{eq:3-42} is bounded in $L^2_{\mathrm{loc}}(\mathbb{R})$, the sequence $w_{i\beta_k}$ is bounded in $H^1_{\mathrm{loc}}(\mathbb{R})$; hence, up to a subsequence, $w_{i\beta_k}\to w_{i0}$ strongly in $L^m_{\mathrm{loc}}(\mathbb{R})$ for every $2\le m<\infty$. Taking the weak limit of \eqref{eq:3-42}, one can find that $(w_{10},w_{20})$ satisfies
\begin{equation}\label{eq:3-46}
	\begin{cases}
		(-\Delta)^{\frac12}w_{10} + w_{10} = a_1 w_{10}^3 + \beta^* w_{20}^2 w_{10},\\
		(-\Delta)^{\frac12}w_{20} + w_{20} = a_2 w_{20}^3 + \beta^* w_{10}^2 w_{20}.
	\end{cases}
\end{equation}
We now show that there exist $r_2,\delta_2>0$ such that
\begin{equation}\label{eq:3-47}
	\liminf_{k\to\infty}\int_{B_{r_2}(0)}|w_{2\beta_k}(x)|^2dx\ge \delta_2>0.
\end{equation}
Indeed, it suffices to prove $w_{20}\neq0$. Arguing by contradiction that $w_{20}=0$, we then see from \eqref{eq:3-46} that $w_{10}\ge0$ satisfies
\begin{equation}\label{eq:3-48}
	(-\Delta)^{\frac12}w_{10} + w_{10} = a_1 w_{10}^3 \quad\text{in }\mathbb R.
\end{equation}
By a simple rescaling, $\sqrt{a_1}\,w_{10}$ satisfies \eqref{singal}; together with \eqref{pohozaev-eq} and \eqref{lem:single_gns} for $p=4$, this yields
\begin{equation}\label{eq:3-49}
\|w_{10}\|_2^2 \ge \frac{a^*}{a_1}>1.
\end{equation}
However, it follows from \eqref{eq:3-19} and weak lower semi-continuity of $\|\cdot\|_2^2$ that
\begin{equation}\label{eq:3-50}
	\|w_{10}\|_2^2 \le \liminf_{k\to\infty}\big(\|w_{1\beta_k}\|_2^2 + \|w_{2\beta_k}\|_2^2\big)=1,
\end{equation}
which leads to an immediate contradiction. Therefore, \eqref{eq:3-47} holds. Combining \eqref{eq:3-35} and \eqref{eq:3-47}, we also obtain \eqref{eq:3-18}, which finishes (i).

(ii) To prove (ii), from \eqref{eq:3-37}, it suffices to show that $x_0\in\mathcal{V}$. Arguing by contradiction that $x_0\notin\mathcal{V}$, without loss of generality, we suppose $V_1(x_0)\neq 0$, which indicates $V_1(x_0)>0$. We then deduce from \eqref{eq:3-18} that
\begin{equation}\label{eq:3-51}
	\liminf_{k\to\infty}\sum_{i=1}^2\int_{\mathbb R} V_i(\varepsilon_{\beta_k}x+x_{\beta_k})|w_{i\beta_k}|^2 dx
	\ge V_1(x_0)\liminf_{k\to\infty}\int_{B_{r_0}(0)}|w_{1\beta_k}|^2 dx\ge V_1(x_0)\delta>0,
\end{equation}
which contradicts \eqref{eq:3-36} so that (ii) follows.

(iii) According to Theorem 1.1 in \cite{RosOtonSerra2014} applied to the whole space $\mathbb R$, the boundary terms vanish. In our specific case $N=1$ and $s=\frac 12$, the coefficient $(2s-N)$ vanishes identically. Thus, we can deduce that
\begin{equation}\label{eq:3-52}
	\int_{\mathbb R}(w_{10}^2+w_{20}^2)dx=\frac12\int_{\mathbb R}\big(a_1w_{10}^4+2\beta^*w_{10}^2w_{20}^2+a_2w_{20}^4\big)dx.
\end{equation}
Linearly combining \eqref{eq:3-52} with the standard Nehari identity of \eqref{eq:3-46}, one directly deduces the following Pohozaev identity
\begin{equation}\label{eq:3-53}
	\begin{aligned}
		&\hphantom{=}\int_{\mathbb R}\big(|(-\Delta)^{\frac14}w_{10}|^2+|(-\Delta)^{\frac14}w_{20}|^2\big)dx
		+\frac12\int_{\mathbb R}(w_{10}^2+w_{20}^2)dx\\
		&=\frac34\int_{\mathbb R}\big(a_1w_{10}^4+2\beta^*w_{10}^2w_{20}^2+a_2w_{20}^4\big)dx.
	\end{aligned}
\end{equation}
Note from \eqref{eq:3-19} that
\begin{equation}\label{eq:3-54}
	\int_{\mathbb R}\big(|(-\Delta)^{\frac14}w_{10}|^2+|(-\Delta)^{\frac14}w_{20}|^2\big)dx\le 1
	\quad \text{and}\quad
	\int_{\mathbb R}(w_{10}^2+w_{20}^2)dx\le 1.
\end{equation}
Moreover, we can infer from \eqref{eq:3-46} that
\begin{equation}\label{eq:3-55}
	\begin{aligned}
		&\hphantom{=}\int_{\mathbb R}\big(|(-\Delta)^{\frac14}w_{10}|^2+|(-\Delta)^{\frac14}w_{20}|^2\big)dx
		+\int_{\mathbb R}(w_{10}^2+w_{20}^2)dx\\
		&=\int_{\mathbb R}\big(a_1w_{10}^4+a_2w_{20}^4+2\beta^*w_{10}^2w_{20}^2\big)dx.
	\end{aligned}
\end{equation}
Combining \eqref{eq:3-53} with \eqref{eq:3-55}, we can infer that
\begin{equation}\label{eq:3-56}
	2\int_{\mathbb R}(w_{10}^2+w_{20}^2)dx=\int_{\mathbb R}\big(a_1w_{10}^4+a_2w_{20}^4+2\beta^*w_{10}^2w_{20}^2\big)dx.
\end{equation}
Together with the refined Gagliardo--Nirenberg inequality \eqref{refine} and Lemma \ref{holder} for $p=4$ and $q=2$, we deduce from \eqref{eq:3-54}--\eqref{eq:3-56} that
\begin{equation}\label{eq:3-57}
\begin{aligned}
    &\hphantom{=}\int_{\mathbb R}\big(|(-\Delta)^{\frac14}w_{10}|^2+|(-\Delta)^{\frac14}w_{20}|^2\big)dx\\
	&=\int_{\mathbb R}(w_{10}^2+w_{20}^2)dx\\
	&=\frac12\int_{\mathbb R}\big(a_1w_{10}^4+a_2w_{20}^4+2\beta^*w_{10}^2w_{20}^2\big)dx\\
	&=\frac{a^*}{2}\int_{\mathbb R}(w_{10}^2+w_{20}^2)^2dx
	-\frac12\Big[\sqrt{a^*-a_1}\sqrt{\int_{\mathbb R}w_{10}^4dx}-\sqrt{a^*-a_2}\sqrt{\int_{\mathbb R}w_{20}^4dx}\Big]^2\\
	&\quad -\sqrt{(a^*-a_1)(a^*-a_2)}\Big[\sqrt{\int_{\mathbb R}w_{10}^4dx}\sqrt{\int_{\mathbb R}w_{20}^4dx}-\int_{\mathbb R}w_{10}^2w_{20}^2dx\Big]\\
    &\le \frac{a^*}{2}\int_{\mathbb R}(w_{10}^2+w_{20}^2)^2dx\\
	&\le \int_{\mathbb R}\big(|(-\Delta)^{\frac14}w_{10}|^2+|(-\Delta)^{\frac14}w_{20}|^2\big)dx\int_{\mathbb R}(w_{10}^2+w_{20}^2)dx\\
	&\le \int_{\mathbb R}\big(|(-\Delta)^{\frac14}w_{10}|^2+|(-\Delta)^{\frac14}w_{20}|^2\big)dx,
	\end{aligned}
\end{equation}
which suggests that the aforementioned equalities hold. Therefore, based on Lemma \ref{holder} and Lemma \ref{lem:optimizers}, we can infer that
\begin{equation}\label{eq:3-58}
	(w_{10}(x),w_{20}(x))=\big(\tau\sin\theta\,Q(\alpha x-\eta_0),\tau\cos\theta\,Q(\alpha x-\eta_0)\big),\ w_{10}^2(x)=k w_{20}^2(x),
\end{equation}
with some $\tau>0$, $\alpha>0$, $\theta\in[0,2\pi)$, $\eta_0\in\mathbb R$ and $k>0$. Actually, it follows from \eqref{eq:3-57} that
\begin{equation}\label{eq:3-59}
	\int_{\mathbb R}\big(|(-\Delta)^{\frac14}w_{10}|^2+|(-\Delta)^{\frac14}w_{20}|^2\big)dx=\int_{\mathbb R}(w_{10}^2+w_{20}^2)dx=1,
\end{equation}
which implies $\alpha=1$ and $\tau=\frac{1}{\|Q\|_{2}}$ in \eqref{eq:3-58}. Besides, \eqref{eq:3-57} also gives
\begin{equation}\label{eq:3-60}
	\sqrt{a^*-a_1}\sqrt{\int_{\mathbb R}w_{10}^4dx}=\sqrt{a^*-a_2}\sqrt{\int_{\mathbb R}w_{20}^4dx},
\end{equation}
which yields $k=\sqrt{\frac{a^*-a_2}{a^*-a_1}}$ in \eqref{eq:3-58} and through simple calculations we thus conclude that
\begin{equation}\label{eq:3-61}
	w_{10}(x)=\frac{\sqrt{\gamma}}{\|Q\|_2}Q(x-\eta_0)
	\quad \text{and}\quad
	w_{20}(x)=\frac{\sqrt{1-\gamma}}{\|Q\|_2}Q(x-\eta_0),
\end{equation}
where $\gamma$ is given by \eqref{eq:3-1}. Moreover, one can further derive from \eqref{eq:3-19} and \eqref{eq:3-59} that
\begin{equation}\label{eq:3-62}
\lim_{k\to\infty}\int_{\mathbb R}\big(|(-\Delta)^{\frac14}w_{i\beta_k}|^2+w_{i\beta_k}^2\big)dx
=\int_{\mathbb R}\big(|(-\Delta)^{\frac14}w_{i0}|^2+w_{i0}^2\big)dx,\ i=1,2.
\end{equation}
Combining with the fact that $w_{i\beta_k}\rightharpoonup w_{i0}$ in $H^{\frac12}(\mathbb R)$, we thus obtain that
\begin{equation}\label{eq:3-63}
	w_{i\beta_k}\to w_{i0} \text{ strongly in } H^{\frac12}(\mathbb R),\ i=1,2.
\end{equation}
By \eqref{eq:3-61}, \eqref{eq:3-24} holds, and the proof of the lemma is complete.
\end{proof}

\begin{proofof}{Theorem \ref{con}}
	Notice that from \eqref{eq:3-42}, one can obtain
\begin{equation}\label{eq:3-64}
	(-\Delta)^{\frac12}w_{i\beta_k}(x)-C_{i\beta_k}(x)w_{i\beta_k}(x)\le 0
	\quad \text{in } \mathbb R,\ i=1,2,
\end{equation}
where $C_{1\beta_k}:=a_1w_{1\beta_k}^2+\beta_k w_{2\beta_k}^2$ and $C_{2\beta_k}:=a_2w_{2\beta_k}^2+\beta_k w_{1\beta_k}^2$. Indeed, the additional term $(\varepsilon_{\beta_k}\mu_{\beta_k}-\varepsilon_{\beta_k}V_i(\varepsilon_{\beta_k}x+x_{\beta_k}))w_{i\beta_k}$ is nonpositive for $k$ large, since $\varepsilon_{\beta_k}\mu_{\beta_k}\to-1$ and $V_i\ge0$.

Thanks to the fractional Sobolev embedding in one dimension, $H^{\frac12}(\mathbb R)\hookrightarrow L^m(\mathbb R)$ for all $m\ge 2$, which guarantees $C_{i\beta_k}\in L^{\corr{r}}(\mathbb R)$ for any $\corr{r}\ge1$. Thus, applying the non-local De Giorgi--Nash--Moser theory (see Theorem 1.1 in \cite{Kim2025}) to the globally nonnegative weak subsolution $w_{i\beta_k}\ (i=1,2)$ of \eqref{eq:3-64},
\begin{equation}\label{eq:3-65}
	\sup_{B_1(\xi_0)} w_{i\beta_k}
	\le C\Big(\int_{B_2(\xi_0)}|w_{i\beta_k}|^2dx\Big)^{\frac12}+\mathrm{Tail}(w_{i\beta_k};\xi_0,1),\ i=1,2,
\end{equation}
where $C$ is a positive constant, $\xi_0$ is an arbitrary point in $\mathbb R$ and
\begin{equation}\label{eq:3-66}
\mathrm{Tail}(w_{i\beta_k};\xi_0,1)=\int_{|x-\xi_0|>1}\frac{w_{i\beta_k}(x)}{|x-\xi_0|^2}dx.
\end{equation}
Moreover, from \eqref{eq:3-24} we can deduce
\begin{equation}\label{eq:3-67}
	\int_{|x|>r}|w_{i\beta_k}|^2dx\to 0
	\quad \text{and}\quad
	\mathrm{Tail}(w_{i\beta_k};\xi_0,\frac r2)\to0 \text{ as } r\to\infty \text{ uniformly in }k,\ i=1,2.
\end{equation}
Combining \eqref{eq:3-65} and \eqref{eq:3-67}, we thus infer that
\begin{equation}\label{eq:3-68}
	w_{i\beta_k}(x)\to0 \quad \text{as } |x|\to\infty \text{ uniformly in }k,\ i=1,2.
\end{equation}
Consequently, each $w_{i\beta_k}\ (i=1,2)$ has at least one global maximum point. Let $z_{i\beta_k}$ be a global maximum point of $u_{i\beta_k}$, $i=1,2$; then, by \eqref{eq:3-17}, $w_{i\beta_k}$ attains its global maximum at the point $\frac{z_{i\beta_k}-x_{\beta_k}}{\varepsilon_{\beta_k}}$. Hence, one can derive from \eqref{eq:3-18} and \eqref{eq:3-68} that
\begin{equation}\label{eq:3-69}
	\limsup_{k\to\infty}\frac{|z_{i\beta_k}-x_{\beta_k}|}{\varepsilon_{\beta_k}}<\infty,\ i=1,2.
\end{equation}
Going if necessary to a subsequence, there exists $\eta_i\in\mathbb R$ such that
\begin{equation}\label{eq:3-70}
	\lim_{k\to\infty}\frac{z_{i\beta_k}-x_{\beta_k}}{\varepsilon_{\beta_k}}=\eta_i,\ i=1,2.
\end{equation}
Define
\begin{equation}\label{eq:3-71}
	\overline{w}_{i\beta_k}(x):=\varepsilon_{\beta_k}^{\frac12}u_{i\beta_k}\big(\varepsilon_{\beta_k}x+z_{i\beta_k}\big)=w_{i\beta_k}\big(x+\frac{z_{i\beta_k}-x_{\beta_k}}{\varepsilon_{\beta_k}}\big),\ i=1,2.
\end{equation}
From Lemma \ref{lemma3-2} (iii), we obtain
\begin{equation}\label{eq:3-72}
\begin{cases}
	\displaystyle\lim_{k\to\infty}\overline{w}_{1\beta_k}(x)=\overline{w}_{10}(x)=\frac{\sqrt{\gamma}}{\|Q\|_2}Q(x+\eta_1-\eta_0),\\[6pt]
	\displaystyle\lim_{k\to\infty}\overline{w}_{2\beta_k}(x)=\overline{w}_{20}(x)=\frac{\sqrt{1-\gamma}}{\|Q\|_2}Q(x+\eta_2-\eta_0),
\end{cases}
\end{equation}
strongly in $H^{\frac12}(\mathbb R)$. Note from \eqref{eq:3-71} that the origin is a critical point of $\overline{w}_{i\beta_k}$; passing to the limit, it is also a critical point of $\overline{w}_{i0}$ for $i=1,2$. Since the ground state $Q$ is radially symmetric and strictly decreasing in $r=|x|$, we deduce that $\overline{w}_{i0}$ attains its maximum at the origin, which forces $\eta_1=\eta_2=\eta_0$ in \eqref{eq:3-72}. Therefore, \eqref{eq:2-18} and \eqref{eq:2-20} immediately follow from \eqref{eq:3-70} and \eqref{eq:3-72}. Furthermore, \eqref{eq:3-68} and Lemma \ref{lemma3-2} (ii) imply that \eqref{eq:2-16} holds. Finally, \eqref{eq:2-21} follows from \eqref{eq:3-20} and \eqref{eq:3-22}, which completes the proof of Theorem \ref{con}.
\end{proofof}

\section{Explicit concentration behavior}
In this section we study the concentration behavior of the minimizers of problem \eqref{eq:min_problem} for $q=2$ and $p=4$ under explicit trapping potentials, namely the polynomial potentials introduced in \eqref{eq:potential_form}; some delicate estimates of the associated functionals will also be given.

\begin{lemma}\label{lem:4-1}
	Suppose that $V_i(x)$ ($i=1,2$) satisfy \eqref{eq:potential_form} and \eqref{eq:common_zeros_condition}. If $0<q_0<1$, then for any given $0<a_1,a_2<a^*$ and $0<\beta<\beta^*$, we have
	\begin{equation}\label{eq:4-1}
		\liminf_{\beta\to\beta^*}\frac{e(a_1,a_2,\beta)}{(\beta^*-\beta)^{\frac{q_0}{q_0+1}}}
		\le \frac{q_0+1}{q_0 a^*}(q_0\lambda_0)^{\frac{1}{q_0+1}}
		\big[2\gamma(1-\gamma)\big]^{\frac{q_0}{q_0+1}},
	\end{equation}
	where $\gamma\in(0,1)$, $q_0>0$ and $\lambda_0>0$ are given in \eqref{eq:3-1}, \eqref{eq:q_j_q_0} and \eqref{eq:lambda_0}, respectively.
\end{lemma}

\begin{proof}
	Let $(\phi_{1\sigma},\phi_{2\sigma})$ be the functions given by \eqref{eq:test_sequence} for $N=1$ with
	\begin{equation}\label{eq:4-2}
		\sigma=\left[\frac{q_0\lambda_0}{2\gamma(1-\gamma)(\beta^*-\beta)}\right]^{\frac{1}{q_0+1}}.
	\end{equation}
	At the same time, set $\theta=\gamma$ and take $x_0=x_{j_0}$, where $x_{j_0}\in Z_0$ with $Z_0$ given by \eqref{eq:set_Z_0}. Applying \eqref{eq:test_energy_final}, one can obtain
	\begin{equation}\label{eq:4-3}
	\begin{aligned}
		&\hphantom{=}\sum_{i=1}^2\int_{\mathbb R}|(-\Delta)^{\frac14}\phi_{i\sigma}|^2dx
		-\frac12\int_{\mathbb R}\big(a_1\phi_{1\sigma}^4+a_2\phi_{2\sigma}^4+2\beta\phi_{1\sigma}^2\phi_{2\sigma}^2\big)dx\\
		&\le \frac{2\sigma}{a^*}\gamma(1-\gamma)(\beta^*-\beta)+O(\sigma^{-1}).
	\end{aligned}
	\end{equation}
	Note that, from $0<q_0<1$, we can deduce that
	\begin{equation}\label{eq:4-4}
		\frac{2\sigma}{a^*}\gamma(1-\gamma)(\beta^*-\beta)+O(\sigma^{-1})
		=\frac{2\sigma}{a^*}\gamma(1-\gamma)(\beta^*-\beta)+o(\sigma^{-q_0}),
	\end{equation}
	as $\sigma\to\infty$. Therefore, combining \eqref{eq:4-3} and \eqref{eq:4-4}, one can infer from \eqref{eq:4-2} that
	\begin{equation}\label{eq:4-5}
		\begin{aligned}
			&\hphantom{=}\sum_{i=1}^2\int_{\mathbb R}|(-\Delta)^{\frac14}\phi_{i\sigma}|^2dx
			-\frac12\int_{\mathbb R}\big(a_1\phi_{1\sigma}^4+a_2\phi_{2\sigma}^4+2\beta\phi_{1\sigma}^2\phi_{2\sigma}^2\big)dx\\
			&\le \frac{1}{ a^*}(q_0\lambda_0)^{\frac{1}{q_0+1}}
			\big[2\gamma(1-\gamma)\big]^{\frac{q_0}{q_0+1}}(\beta^*-\beta)^{\frac{q_0}{q_0+1}}
			+o\big((\beta^*-\beta)^{\frac{q_0}{q_0+1}}\big) \text{ as } \beta\to\beta^*.
		\end{aligned}
	\end{equation}
	On the other hand, without loss of generality, we suppose $q_0=q_{1j_0}<q_{2j_0}$ (the case $q_0=q_{2j_0}<q_{1j_0}$ is symmetric; if $q_{1j_0}=q_{2j_0}=q_0$, both terms below are of order $\sigma^{-q_0}$ and the second branch of \eqref{eq:lambda_j} gives the same conclusion), hence we can deduce from \eqref{eq:test_energy_final} and \eqref{eq:lambda_j}--\eqref{eq:set_Z_0} that
	\begin{equation}\label{eq:4-6}
		\begin{aligned}
			&\hphantom{=}\int_{\mathbb R}V_1(x)\phi_{1\sigma}^2dx+\int_{\mathbb R}V_2(x)\phi_{2\sigma}^2dx\\
			&=\sigma^{-q_0}\int_{\mathbb R}\frac{A_\sigma^2\gamma V_1(\sigma^{-1}x+x_{j_0})}{a^*}\sigma^{q_0}\varphi^2(\sigma^{-1}x)Q^2(x)dx\\
			&\quad+\sigma^{-q_0}\int_{\mathbb R}\frac{A_\sigma^2(1-\gamma)V_2(\sigma^{-1}x+x_{j_0})}{a^*}\sigma^{q_0}\varphi^2(\sigma^{-1}x)Q^2(x)dx\\
			&=\sigma^{-q_0}\int_{\mathbb R}\frac{A_\sigma^2\gamma V_1(\sigma^{-1}x+x_{j_0})}{a^*|\sigma^{-1}x|^{q_{1j_0}}}|x|^{q_{1j_0}}\varphi^2(\sigma^{-1}x)Q^2(x)dx\\
			&\quad+\sigma^{-q_0}\int_{\mathbb R}\frac{A_\sigma^2(1-\gamma) V_2(\sigma^{-1}x+x_{j_0})}{a^*|\sigma^{-1}x|^{q_{2j_0}}}\sigma^{q_{0}-q_{2j_0}}|x|^{q_{2j_0}}\varphi^2(\sigma^{-1}x)Q^2(x)dx\\
			&=\frac{\lambda_0}{a^*\sigma^{q_0}}+o(\sigma^{-q_0}) \quad\text{as } \sigma\to\infty.
		\end{aligned}
	\end{equation}
	Therefore, from \eqref{eq:4-2}, we thus infer
	\begin{equation}\label{eq:4-7}
		\begin{aligned}
			&\hphantom{=}\int_{\mathbb R}V_1(x)\phi_{1\sigma}^2dx+\int_{\mathbb R}V_2(x)\phi_{2\sigma}^2dx\\
			&=\frac{1}{q_0 a^*}(q_0\lambda_0)^{\frac{1}{q_0+1}}
			\big[2\gamma(1-\gamma)\big]^{\frac{q_0}{q_0+1}}(\beta^*-\beta)^{\frac{q_0}{q_0+1}}
			+o\big((\beta^*-\beta)^{\frac{q_0}{q_0+1}}\big)\quad\text{as } \beta\to\beta^*.
		\end{aligned}
	\end{equation}
	Combining \eqref{eq:4-5} and \eqref{eq:4-7}, we then derive the estimate \eqref{eq:4-1}, as claimed.
\end{proof}

\begin{proofof}{Theorem \ref{thm:1.3}}
	Let $\{(u_{1\beta_k},u_{2\beta_k})\}$ be the subsequence obtained in Theorem \ref{con}. Applying the refined Gagliardo--Nirenberg inequality \eqref{refine} and H\"older's inequality to \eqref{eq:3-9}, we then obtain for any $k\in\mathbb N$
\begin{equation}\label{eq:4-8}
	\begin{aligned}
		e(a_1,a_2,\beta_k)
		&\ge (\beta^*-\beta_k)\int_{\mathbb R}u_{1\beta_k}^2u_{2\beta_k}^2dx
		+\sum_{i=1}^2\int_{\mathbb R}V_i(x)u_{i\beta_k}^2dx\\
		&=\frac{(\beta^*-\beta_k)}{\varepsilon_{\beta_k}}\int_{\mathbb R}\overline{w}_{1\beta_k}^2\overline{w}_{2\beta_k}^2dx
		+\sum_{i=1}^2\int_{\mathbb R}V_i(\varepsilon_{\beta_k}x+z_{i\beta_k})\overline{w}_{i\beta_k}^2dx,
	\end{aligned}
\end{equation}
where the second equality comes from \eqref{eq:3-71}. Using \eqref{eq:3-22} and \eqref{eq:3-71}, we find that
\begin{equation}\label{eq:4-9}
	\int_{\mathbb R}\overline{w}_{1\beta_k}^2\overline{w}_{2\beta_k}^2dx
	=\frac{2\gamma(1-\gamma)}{a^*}+o(1) \quad\text{as }k\to\infty.
\end{equation}
Note that $x_0\in \mathcal{V}$ in \eqref{eq:2-16}; hence there exists $1\le j_0\le b$ such that $x_0=x_{j_0}$. Without loss of generality, we may assume that $q_{j_0}=q_{1j_0}<q_{2j_0}$. Next, we claim that
\begin{equation}\label{eq:4-10}
	\frac{z_{i\beta_k}-x_{j_0}}{\varepsilon_{\beta_k}} \text{ is bounded as }k\to\infty,
\end{equation}
and
\begin{equation}\label{eq:4-11}
	q_{j_0}=q_0.
\end{equation}

Suppose by contradiction that either $q_{j_0}<q_0$ or \eqref{eq:4-10} is false. Then, for any given $M>0$, one can infer from Fatou's Lemma that
\begin{equation}\label{eq:4-12}
	\begin{aligned}
	&\hphantom{=}\liminf_{k\to\infty}\varepsilon_{\beta_k}^{-q_0}\int_{\mathbb R}V_1(\varepsilon_{\beta_k}x+z_{1\beta_k})\overline{w}_{1\beta_k}^2dx\\
		&\ge \int_{\mathbb R}\liminf_{k\to\infty}\varepsilon_{\beta_k}^{q_{j_0}-q_0}
		\left(
		\frac{V_1(\varepsilon_{\beta_k}x+z_{1\beta_k})}{|\varepsilon_{\beta_k}x+z_{1\beta_k}-x_{j_0}|^{q_{j_0}}}
		\left|x+\frac{z_{1\beta_k}-x_{j_0}}{\varepsilon_{\beta_k}}\right|^{q_{j_0}}
		\overline{w}_{1\beta_k}^2
		\right)dx
		\ge M.
	\end{aligned}
\end{equation}
From \eqref{eq:4-8}, \eqref{eq:4-9} and \eqref{eq:4-12}, we thus deduce that
\begin{equation}\label{eq:4-13}
	e(a_1,a_2,\beta_k)
	\ge \frac{2\gamma(1-\gamma)}{a^*}(\beta^*-\beta_k)\varepsilon_{\beta_k}^{-1}+M\varepsilon_{\beta_k}^{q_0}
	\ge C M^{\frac{1}{q_0+1}} (\beta^*-\beta_k)^{\frac{q_0}{q_0+1}} \text{ for large } k.
\end{equation}
In fact, the second inequality above holds due to Young's inequality, where $C$ is a suitable positive constant. Since $M$ is arbitrary, the claim follows from the contradiction between \eqref{eq:4-1} and \eqref{eq:4-13}. Following \eqref{eq:2-20} and \eqref{eq:4-10}, after passing to a subsequence, there exists $z_0\in\mathbb R$ such that
\begin{equation}\label{eq:4-14}
	\lim_{k\to\infty}\frac{z_{i\beta_k}-x_{j_0}}{\varepsilon_{\beta_k}}=z_0 \quad\text{for }i=1,2.
\end{equation}

Combining this with \eqref{eq:3-72}, we have
\begin{equation}\label{eq:4-15}
	\begin{aligned}
		&\lim_{k\to\infty}\varepsilon_{\beta_k}^{-q_0}\left(
		\int_{\mathbb R}V_1(\varepsilon_{\beta_k}x+z_{1\beta_k})\overline{w}_{1\beta_k}^2dx
		+\int_{\mathbb R}V_2(\varepsilon_{\beta_k}x+z_{2\beta_k})\overline{w}_{2\beta_k}^2dx
		\right)\\
		=&\lim_{k\to\infty}\int_{\mathbb R}\frac{V_1(\varepsilon_{\beta_k}x+z_{1\beta_k})}{|\varepsilon_{\beta_k}x+z_{1\beta_k}-x_{j_0}|^{q_{1j_0}}}
		\left|x+\frac{z_{1\beta_k}-x_{j_0}}{\varepsilon_{\beta_k}}\right|^{q_{1j_0}}\overline{w}_{1\beta_k}^2dx\\
		&+\lim_{k\to\infty}\varepsilon_{\beta_k}^{q_{2j_0}-q_0}\int_{\mathbb R}\frac{V_2(\varepsilon_{\beta_k}x+z_{2\beta_k})}{|\varepsilon_{\beta_k}x+z_{2\beta_k}-x_{j_0}|^{q_{2j_0}}}
		\left|x+\frac{z_{2\beta_k}-x_{j_0}}{\varepsilon_{\beta_k}}\right|^{q_{2j_0}}\overline{w}_{2\beta_k}^2dx\\
		=&\lim_{x\to x_{j_0}}\frac{\gamma V_1(x)}{a^*|x-x_{j_0}|^{q_{1j_0}}}\int_{\mathbb R}|y+z_0|^{q_{1j_0}}Q^2(y)dy.
	\end{aligned}
\end{equation}
Note that
\begin{equation}\label{eq:4-16}
	\int_{\mathbb R}|y+z_0|^{q_{1j_0}}Q^2(y)dy
	\ge \int_{\mathbb R}|y|^{q_{1j_0}}Q^2(y)dy,
\end{equation}
where the equality holds if and only if $z_0=0$. Indeed, since $Q$ is radially symmetric and decreasing, we have
\begin{equation}\label{eq:4-17}
	\int_{\mathbb R}|y+z_0|^{q_{1j_0}}Q^2(y)dy
	=\int_{\mathbb R}|y-z_0|^{q_{1j_0}}Q^2(y)dy,
\end{equation}
then, we obtain
\begin{equation}\label{eq:4-18}
	\begin{aligned}
		&\int_{\mathbb R}|y+z_0|^{q_{1j_0}}Q^2(y)dy-\int_{\mathbb R}|y|^{q_{1j_0}}Q^2(y)dy\\
		=&\frac12\int_{\mathbb R}\big(|y+z_0|^{q_{1j_0}}-|y|^{q_{1j_0}}\big)Q^2(y)dy
		+\frac12\int_{\mathbb R}\big(|y-z_0|^{q_{1j_0}}-|y|^{q_{1j_0}}\big)Q^2(y)dy\\
		=&\frac12\int_{\mathbb R}\big(|y+z_0|^{q_{1j_0}}-|y|^{q_{1j_0}}\big)Q^2(y)dy
		+\frac12\int_{\mathbb R}\big(|y|^{q_{1j_0}}-|y+z_0|^{q_{1j_0}}\big)Q^2(y+z_0)dy\\
		=&\frac12\int_{\mathbb R}\big(|y+z_0|^{q_{1j_0}}-|y|^{q_{1j_0}}\big)\big(Q^2(y)-Q^2(y+z_0)\big)dy
		\ge 0,
	\end{aligned}
\end{equation}
which implies that \eqref{eq:4-16} holds. Therefore, from the definition of $\lambda_{j_0}$ and $\lambda_0$, one can find that
\begin{equation}\label{eq:4-19}
	\begin{aligned}
		\lim_{x\to x_{j_0}}\frac{\gamma V_1(x)}{a^*|x-x_{j_0}|^{q_{1j_0}}}\int_{\mathbb R}|y+z_0|^{q_{1j_0}}Q^2(y)dy
		&\ge \lim_{x\to x_{j_0}}\frac{\gamma V_1(x)}{a^*|x-x_{j_0}|^{q_{1j_0}}}\int_{\mathbb R}|y|^{q_{1j_0}}Q^2(y)dy\\
		&=\frac{\lambda_{j_0}}{a^*}\ge\frac{\lambda_0}{a^*},
	\end{aligned}
\end{equation}
where the equalities hold if and only if $z_0=0$ and $x_{j_0}\in Z_0$. Applying \eqref{eq:4-8}, \eqref{eq:4-9}, \eqref{eq:4-15} and \eqref{eq:4-19}, we further obtain
\begin{equation}\label{eq:4-20}
	e(a_1,a_2,\beta_k)
	\ge \left(\frac{2\gamma(1-\gamma)}{a^*}+o(1)\right)(\beta^*-\beta_k)\varepsilon_{\beta_k}^{-1}
	+\left(\frac{\lambda_0}{a^*}+o(1)\right)\varepsilon_{\beta_k}^{q_0}
	\quad\text{as }k\to\infty.
\end{equation}
Therefore, one can verify by taking the infimum of the above inequality that
\begin{equation}\label{eq:4-21}
	\liminf_{k\to\infty}\frac{e(a_1,a_2,\beta_k)}{(\beta^*-\beta_k)^{\frac{q_0}{q_0+1}}}
	\ge \frac{q_0+1}{q_0 a^*}(q_0\lambda_0)^{\frac{1}{q_0+1}}
	\big[2\gamma(1-\gamma)\big]^{\frac{q_0}{q_0+1}},
\end{equation}
where the equality holds if and only if $z_0=0$, $x_{j_0}\in Z_0$ and
\begin{equation}\label{eq:4-22}
	\lim_{k\to\infty}\frac{\varepsilon_{\beta_k}}{\epsilon_{\beta_k}}=1 \text{ with } \epsilon_{\beta_k}>0 \text{ given by \eqref{eq:epsilon_beta_k}}.
\end{equation}
Note from \eqref{eq:4-1} and \eqref{eq:4-21} that the equality in \eqref{eq:4-21} holds, hence $z_0=0$, $x_{j_0}\in Z_0$ and \eqref{eq:4-22} holds. Consequently, \eqref{eq:profile_limit} and \eqref{eq:energy_limit_1} follow from \eqref{eq:2-18}, \eqref{epsilo} and \eqref{eq:4-22}. Moreover, using the fact that $z_0=0$, $x_{j_0}\in Z_0$, \eqref{eq:z_limit} follows from \eqref{eq:4-14} and \eqref{eq:4-22}. This completes the proof of Theorem \ref{thm:1.3}.
\end{proofof}

\appendix

\section{Equivalence between ground states and constraint minimizers}

In this appendix, we establish the equivalence between ground states of the coupled fractional Schr\"{o}dinger system \eqref{eq:sys} and the constraint minimizers of the variational problem \eqref{eq:min_problem}.

For any given $(a_1, a_2, \beta) \in \mathbb{R}^+ \times \mathbb{R}^+ \times \mathbb{R}^+$, the action functional of system \eqref{eq:sys} associated with $\mu \in \mathbb{R}$ is defined by
\begin{equation}\label{eq:A.1}
	\begin{aligned}
		I_\mu(u_1, u_2) &:= \frac{1}{2} \sum_{i=1}^2 \left( \|(-\Delta)^{\frac s2} u_i\|_2^2 + \int_{\mathbb{R}^N} (V_i(x) - \mu) u_i^2 \, dx \right) \\
		&\quad - \frac{1}{p} \int_{\mathbb{R}^N} \left( a_1 |u_1|^p + a_2 |u_2|^p + 2\beta |u_1|^q |u_2|^q \right) dx \\
		&= \frac{1}{2} E_{a_1, a_2, \beta}(u_1, u_2) - \frac{\mu}{2} \int_{\mathbb{R}^N} (u_1^2 + u_2^2) \, dx, \quad (u_1, u_2) \in \mathcal{H}_1 \times \mathcal{H}_2,
	\end{aligned}
\end{equation}
where $p = 2 + \frac{4s}{N}$, $q = \frac{p}{2} = 1 + \frac{2s}{N}$, and $E_{a_1, a_2, \beta}(u_1, u_2)$ is given by \eqref{eq:energy}.

The set of all nontrivial weak solutions for system \eqref{eq:sys} is given by
\begin{equation}\label{eq:A.2}
	\mathcal{S}_\mu := \left\{ (u_1, u_2) \in (\mathcal{H}_1 \times \mathcal{H}_2) \setminus\{(0,0)\} : \langle I_\mu'(u_1, u_2), (\varphi_1, \varphi_2) \rangle = 0,\ \forall (\varphi_1, \varphi_2) \in \mathcal{H}_1 \times \mathcal{H}_2 \right\},
\end{equation}
and the set of all ground states for system \eqref{eq:sys} is defined as
\begin{equation}\label{eq:A.3}
	\mathcal{G}_\mu := \left\{ (u_1, u_2) \in \mathcal{S}_\mu : I_\mu(u_1, u_2) \le I_\mu(\bar{u}_1, \bar{u}_2) \text{ for all } (\bar{u}_1, \bar{u}_2) \in \mathcal{S}_\mu \right\}.
\end{equation}

\begin{proposition}\label{Proposition A.1.}
	Suppose $(a_1,a_2,\beta)\in\mathbb{R}^+\times\mathbb{R}^+\times\mathbb{R}^+$ is given, and let $(u_{1\beta},u_{2\beta})$ be a minimizer of problem \eqref{eq:min_problem} with Lagrange multiplier $\mu\in\mathbb{R}$. Then every ground state in $\mathcal{G}_\mu$ is a minimizer of problem \eqref{eq:min_problem}; conversely, $(u_{1\beta},u_{2\beta})\in\mathcal{G}_\mu$.
\end{proposition}

\begin{proof}
	For any given $(a_1, a_2, \beta) \in \mathbb{R}^+ \times \mathbb{R}^+ \times \mathbb{R}^+$, assume that $(u_{1\beta}, u_{2\beta})$ is a minimizer of \eqref{eq:min_problem}. By the Lagrange multiplier theorem, there exists $\mu \in \mathbb{R}$ such that $(u_{1\beta}, u_{2\beta}) \in \mathcal{S}_\mu$; fix this $\mu$. Let $(u_1, u_2)$ be an arbitrary ground state in $\mathcal{G}_\mu$. Set
	\begin{equation}\label{eq:A.4}
		\tilde{u}_i(x) := \frac{u_i(x)}{\sqrt{\rho}}, \quad \text{where } \rho := \int_{\mathbb{R}^N} (u_1^2 + u_2^2) \, dx > 0, \quad i = 1, 2,
	\end{equation}
	so that $\int_{\mathbb{R}^N} (\tilde{u}_1^2 + \tilde{u}_2^2) \, dx = 1$, which implies $(\tilde{u}_1, \tilde{u}_2) \in \mathcal{M}$. By \eqref{eq:min_problem} and \eqref{eq:A.3}, we obtain
	\begin{equation}\label{eq:A.5}
		E_{a_1, a_2, \beta}(\tilde{u}_1, \tilde{u}_2) \ge E_{a_1, a_2, \beta}(u_{1\beta}, u_{2\beta}) \quad \text{and} \quad I_\mu(u_{1\beta}, u_{2\beta}) \ge I_\mu(u_1, u_2).
	\end{equation}
	Combining \eqref{eq:A.1}, \eqref{eq:A.4} and \eqref{eq:A.5} yields
	\begin{equation}\label{eq:A.6}
		I_\mu(\tilde{u}_1, \tilde{u}_2) \ge I_\mu(u_{1\beta}, u_{2\beta}) \ge I_\mu(u_1, u_2).
	\end{equation}
	On the other hand, one can infer from \eqref{eq:sys} and \eqref{eq:A.1} that
	\begin{equation}\label{eq:A.7}
		I_\mu(u_1, u_2) = \left( \frac{1}{2} - \frac{1}{p} \right) \int_{\mathbb{R}^N} \left( a_1 |u_1|^p + a_2 |u_2|^p + 2\beta |u_1|^q |u_2|^q \right) dx.
	\end{equation}
	Combining \eqref{eq:sys} and \eqref{eq:A.4}, we also have
	\begin{equation}\label{eq:A.8}
		I_\mu(\tilde{u}_1, \tilde{u}_2) = \left( \frac{1}{2\rho} - \frac{1}{p\rho^{p/2}} \right) \int_{\mathbb{R}^N} \left( a_1 |u_1|^p + a_2 |u_2|^p + 2\beta |u_1|^q |u_2|^q \right) dx.
	\end{equation}
	Consider the scalar function $h(\rho) := \frac{1}{2\rho} - \frac{1}{p\rho^{p/2}}$ for $\rho > 0$. Differentiating $h(\rho)$ yields
	\begin{equation*}
		h'(\rho) = -\frac{1}{2\rho^2} + \frac{1}{2\rho^{p/2+1}} = \frac{1}{2\rho^2} \left( \rho^{1 - p/2} - 1 \right).
	\end{equation*}
	Since $p = 2 + \frac{4s}{N} > 2$, the exponent $1 - p/2 = -\frac{2s}{N} < 0$. Thus $h'(\rho) > 0$ for $\rho \in (0, 1)$ and $h'(\rho) < 0$ for $\rho \in (1, \infty)$. Therefore, $h(\rho)$ attains its unique global maximum on $(0, \infty)$ strictly at $\rho = 1$, giving
	\begin{equation}\label{eq:A.9}
		h(\rho) \le h(1) = \frac{1}{2} - \frac{1}{p},
	\end{equation}
	with equality holding if and only if $\rho = 1$.
	Together with \eqref{eq:A.7}, \eqref{eq:A.8} and \eqref{eq:A.9}, we deduce that
	\begin{equation}\label{eq:A10}
		I_\mu(\tilde{u}_1, \tilde{u}_2) \le I_\mu(u_1, u_2).
	\end{equation}
	One can easily see from \eqref{eq:A.6} and \eqref{eq:A10} that $I_\mu(\tilde{u}_1, \tilde{u}_2) = I_\mu(u_1, u_2)$, and thus $\rho = 1$.
	This implies that all inequalities in \eqref{eq:A.5} and \eqref{eq:A.6} are equalities, i.e.,
	\begin{equation}\label{eq:A.11}
		I_\mu(u_1, u_2) = I_\mu(u_{1\beta}, u_{2\beta}) \quad \text{and} \quad E_{a_1, a_2, \beta}(u_1, u_2) = E_{a_1, a_2, \beta}(u_{1\beta}, u_{2\beta}).
	\end{equation}
	Consequently, $(u_1,u_2)\in\mathcal{M}$ and, by \eqref{eq:A.11},
	\[
		E_{a_1,a_2,\beta}(u_1,u_2)=E_{a_1,a_2,\beta}(u_{1\beta},u_{2\beta})=e(a_1,a_2,\beta),
	\]
	that is, $(u_1,u_2)$ is a constraint minimizer of problem \eqref{eq:min_problem}. Since $(u_1,u_2)\in\mathcal{G}_\mu$ was arbitrary, every ground state in $\mathcal{G}_\mu$ is a minimizer of \eqref{eq:min_problem}.
It remains to prove the converse: $(u_{1\beta},u_{2\beta})\in\mathcal{G}_\mu$. For $(v_1,v_2)\in\mathcal{S}_\mu$, set $\rho_v:=\int_{\mathbb{R}^N}(v_1^2+v_2^2)\,dx$ and $\tilde v_i:=v_i/\sqrt{\rho_v}$. Since $(\tilde v_1,\tilde v_2)\in\mathcal{M}$, \eqref{eq:min_problem} and \eqref{eq:A.1} yield $I_\mu(\tilde v_1,\tilde v_2)\ge I_\mu(u_{1\beta},u_{2\beta})$; on the other hand, \eqref{eq:A.7}--\eqref{eq:A.9} with $(v,\rho_v)$ in place of $(u,\rho)$ yield $I_\mu(\tilde v_1,\tilde v_2)\le I_\mu(v_1,v_2)$. Hence $I_\mu(u_{1\beta},u_{2\beta})\le I_\mu(v_1,v_2)$ for every $(v_1,v_2)\in\mathcal{S}_\mu$, that is, $(u_{1\beta},u_{2\beta})\in\mathcal{G}_\mu$. The proof is complete.

\end{proof}

\section{Convexity inequality for the fractional Laplacian}

\begin{lemma} \label{lem:convexity}
For any $u_1,u_2\in\Hd$, let $v=\sqrt{u_1^2+u_2^2}$. Then $v\in\Hd$ and
\begin{equation}\label{convex}
\norm{(-\Delta)^{\frac s2} v}{2}^2 \leq \norm{(-\Delta)^{\frac s2} u_1}{2}^2 + \norm{(-\Delta)^{\frac s2} u_2}{2}^2.
\end{equation}
\end{lemma}

\begin{proof}
For any $x,y\in\R^N$, set $a=(u_1(x),u_2(x))$, $b=(u_1(y),u_2(y))$. Then $v(x)=\abs{a}$, $v(y)=\abs{b}$. By the reverse triangle inequality for vectors,
\begin{equation}
\abs{v(x)-v(y)} = \abs{ \abs{a}-\abs{b} } \leq \abs{a-b} = \sqrt{ \abs{u_1(x)-u_1(y)}^2 + \abs{u_2(x)-u_2(y)}^2 }.
\end{equation}
Squaring both sides,
\begin{equation}
\abs{v(x)-v(y)}^2 \leq \abs{u_1(x)-u_1(y)}^2 + \abs{u_2(x)-u_2(y)}^2.
\end{equation}
Dividing by $\abs{x-y}^{N+2s}$ and integrating over $\R^{2N}$, we obtain
\begin{equation}
\int_{\R^{2N}} \frac{\abs{v(x)-v(y)}^2}{\abs{x-y}^{N+2s}} dx dy \leq \int_{\R^{2N}} \frac{\abs{u_1(x)-u_1(y)}^2}{\abs{x-y}^{N+2s}} dx dy + \int_{\R^{2N}} \frac{\abs{u_2(x)-u_2(y)}^2}{\abs{x-y}^{N+2s}} dx dy.
\end{equation}
Multiplying by the normalization constant $\frac{C_{N,s}}{2}$ and using Proposition~3.6 in \cite{DiNezzaPalatucciValdinoci2012}, which gives
\begin{equation}
[u]_{H^s(\mathbb R^N)}^2
=2C_{N,s}^{-1}\|(-\Delta)^{\frac s2}u\|_2^2,
\end{equation}
we obtain the desired inequality. Moreover, \(v\in L^2\), since \(v^2=u_1^2+u_2^2\) and \(u_i\in L^2\) for \(i=1,2\). Since the Gagliardo seminorm of \(v\) is finite by the preceding inequality, we conclude that \(v\in\Hd\).
\end{proof}

\section{Pohozaev identity for the ground state $Q$}

\begin{lemma} \label{pohozaev}
Let $Q$ be the unique positive radial ground state of the single equation \eqref{singal}. Then
\begin{equation}\label{pohozaev-eq}
\norm{(-\Delta)^{\frac s2} Q}{2}^2 = \norm{Q}{2}^2, \quad \norm{Q}{p}^p = \frac{N+2s}{N} \norm{Q}{2}^2.
\end{equation}
\end{lemma}

\begin{proof}
We derive the Pohozaev identity for a solution of the single equation $\fl Q + \frac{2s}{N} Q - \abs{Q}^{\frac{4s}{N}} Q = 0$. Multiply the equation by $x\cdot\nabla Q$ and integrate over $\R^N$. Using the identity (see e.g. \cite{Secchi2016})
\begin{equation} \label{eq:pohozaev_fl}
\intR \fl Q (x\cdot\nabla Q) dx = \frac{2s-N}{2} \norm{(-\Delta)^{\frac s2} Q}{2}^2,
\end{equation}
and
\begin{equation} \label{eq:pohozaev_linear}
\intR \frac{2s}{N} Q (x\cdot\nabla Q) dx = -\frac{2s}{N} \cdot \frac{N}{2} \norm{Q}{2}^2 = -s \norm{Q}{2}^2,
\end{equation}
\begin{equation} \label{eq:pohozaev_nonlinear}
\intR \abs{Q}^{\frac{4s}{N}} Q (x\cdot\nabla Q) dx = -\frac{N}{2+\frac{4s}{N}} \norm{Q}{p}^p,
\end{equation}
we obtain
\begin{equation} \label{eq:pohozaev_first}
\frac{2s-N}{2} \norm{(-\Delta)^{\frac s2} Q}{2}^2 - s \norm{Q}{2}^2 + \frac{N}{2+\frac{4s}{N}} \norm{Q}{p}^p = 0.
\end{equation}

On the other hand, multiplying the original equation by $Q$ and integrating gives
\begin{equation} \label{eq:pohozaev_second}
\norm{(-\Delta)^{\frac s2} Q}{2}^2 + \frac{2s}{N} \norm{Q}{2}^2 = \norm{Q}{p}^p.
\end{equation}

Eliminating $\norm{Q}{p}^p$ from these two relations yields
\begin{equation} \label{eq:pohozaev_kinetic_mass}
\norm{(-\Delta)^{\frac s2} Q}{2}^2 = \norm{Q}{2}^2,
\end{equation}
then, we have
\begin{equation} \label{eq:pohozaev_kinetic_nonlinear}
\norm{(-\Delta)^{\frac s2} Q}{2}^2 = \frac{N}{N+2s} \norm{Q}{p}^p,
\end{equation}
which is equivalent to
\begin{equation}
\norm{Q}{p}^p = \frac{N+2s}{N} \norm{Q}{2}^2.
\end{equation}
The integration by parts for the fractional Laplacian used above can be justified by approximation and the decay of $Q$. This completes the proof.
\end{proof}

\section{Estimate for the truncated ground state}
\begin{lemma}\label{app:trunc-N1}
Let $Q$ be the ground state of \eqref{singal} and let $\varphi_\sigma(x) = \varphi(\sigma^{-1}x)$ be the cut-off function defined in \eqref{varphi}. Then for sufficiently large $\sigma$, it holds that
\begin{equation}\label{eq:A4_estimate}
	\int_{\mathbb{R}^N} |(-\Delta)^{\frac s2}(\varphi_\sigma Q)|^2 dx \le \int_{\mathbb{R}^N} |(-\Delta)^{\frac s2}Q|^2 dx + O(\sigma^{-4s}).
\end{equation}
\end{lemma}
\begin{proof}
We first note that for dimensions $N \ge 2$, the desired estimate \eqref{eq:A4_estimate} has been established in \cite[Lemma 3.2]{DuTianWangZhang2019}. Therefore, it suffices to consider the one-dimensional case $N=1$.

To this end, let \(\varphi_\sigma(x)=\varphi(\sigma ^{-1}x)\), where \(\varphi\in C_0^\infty(\mathbb{R})\) satisfies \eqref{varphi}. Set $Q_\sigma(x):=\varphi_\sigma(x)Q(x).$
Denote
\begin{equation}
J(u):=\int_{\mathbb{R}^{2}}\frac{|u(x)-u(y)|^2}{|x-y|^{1+2s}}\,dx\,dy.
\end{equation}
We use the decomposition in the proof of \cite[Lemma~3.2]{DuTianWangZhang2019}:
\begin{equation}\label{total}
J(Q_\sigma)
=
\int_{B_\sigma\times B_\sigma}
+2\int_{\mathbb{D}}
+2\int_{\mathbb{E}}
+\int_{B_\sigma^c\times B_\sigma^c},
\end{equation}
where
\begin{equation}
\mathbb{D}=\{(x,y):x\in B_\sigma,\ y\in B_\sigma^c,\ |x-y|>\sigma/2\},
\end{equation}
\begin{equation}
\mathbb{E}=\{(x,y):x\in B_\sigma,\ y\in B_\sigma^c,\ |x-y|\le\sigma/2\}.
\end{equation}
The integrals over \(B_\sigma\times B_\sigma\), \(\mathbb{D}\), and \(B_\sigma^c\times B_\sigma^c\) are estimated exactly as in \cite[Lemma~3.2]{DuTianWangZhang2019}; their total contribution is at most
\begin{equation}\label{interval1}
J(Q)+O(\sigma^{-4s}).
\end{equation}
Hence it remains to estimate the integral over \(\mathbb{E}\).

For \((x,y)\in\mathbb{E}\), we have \(y\in B_\sigma^c\) and \(|x-y|\le\sigma/2\), so the segment joining \(x\) and \(y\) remains at distance at least $c\sigma$ from the origin. By \cite[Theorem~1.1]{DuTianWangZhang2019} and the regularity and decay estimates for ground states of fractional Schrödinger equations \cite{FrankLenzmann2013,FrankLenzmannSilvestre2016}, one has
\begin{equation}
|\nabla Q(z)|\le C(1+|z|)^{-2-2s}.
\end{equation}
Since \(\varphi\in C_0^\infty(\mathbb{R})\), there exists a constant \(C>0\), independent of \(\sigma\), such that \(|\varphi'|\le C\). Hence
\begin{equation}
|\nabla\varphi_\sigma(z)|\le C\sigma^{-1}.
\end{equation}
Consequently, for \(|z|\ge\sigma/2\),
\begin{equation}
|\nabla Q_\sigma(z)|
\le |\nabla\varphi_\sigma(z)|Q(z)+\varphi_\sigma(z)|\nabla Q(z)|
\le C\sigma^{-2-2s}.
\end{equation}
Therefore, for \((x,y)\in\mathbb{E}\),
\begin{equation}
|Q_\sigma(x)-Q_\sigma(y)|
\le C\sigma^{-2-2s}|x-y|.
\end{equation}
It follows that
\begin{equation}\label{interval2}
\begin{aligned}
	\int_{\mathbb{E}}\frac{|Q_\sigma(x)-Q_\sigma(y)|^2}{|x-y|^{1+2s}}\,dx\,dy
	&\le C\sigma^{-4-4s}
	\int_{|x|\le\sigma}dx
	\int_{|\xi|\le\sigma/2}|\xi|^{1-2s}\,d\xi \\
	&=O(\sigma^{-1-6s})
	=O(\sigma^{-4s}),
\end{aligned}
\end{equation}
where we used \(s\in(0,1)\) and \(\sigma\to\infty\). Combining \eqref{total},\eqref{interval1} and \eqref{interval2}, one can infer that
\begin{equation}
J(Q_\sigma)\le J(Q)+O(\sigma^{-4s}).
\end{equation}
Finally, since
\begin{equation}
\int_{\mathbb{R}}|(-\Delta)^{\frac s2}u|^2\,dx
=
c_{1,s}J(u)
\end{equation}
for a positive constant \(c_{1,s}\) depending only on \(s\) (see Proposition~3.6 in \cite{DiNezzaPalatucciValdinoci2012}), we obtain
\begin{equation}
\int_{\mathbb{R}}\big|(-\Delta)^{\frac s2}(\varphi_\sigma Q)\big|^2\,dx
\le
\int_{\mathbb{R}}\big|(-\Delta)^{\frac s2}Q\big|^2\,dx
+O(\sigma^{-4s}).
\end{equation}
This completes the proof for \(N=1\).
\end{proof}
\section{Optimizers of the coupled GNS inequality}

\begin{lemma} \label{lem:optimizers}
Let $(u_1,u_2)\in\Hd\times\Hd\setminus\{(0,0)\}$ be such that equality holds in both \eqref{eq:holder_coupled} and \corr{\eqref{eq:nonlinear}}. Then there exist constants $\kappa>0$, $\alpha>0$, $\eta\in\R^N$, and $\theta\in[0,2\pi)$ such that
\begin{equation}
(u_1(x),u_2(x)) = (\kappa \sin\theta \, Q(\alpha x+\eta), \kappa \cos\theta \, Q(\alpha x+\eta)),
\end{equation}
where $Q$ is the ground state of \eqref{singal}. Moreover, if the pair additionally satisfies $\int (u_1^2+u_2^2) dx=1$, then $\kappa^2 \alpha^{-N} \norm{Q}{2}^2=1$.
\end{lemma}

\begin{proof}
\textit{Step 1: Equality in \eqref{eq:holder_coupled} implies proportionality.} By Hölder's inequality, equality holds in \eqref{eq:holder_coupled} if and only if $\abs{u_1}^p$ and $\abs{u_2}^p$ are linearly dependent, i.e., there exists $c\geq0$ such that $\abs{u_2}=c\abs{u_1}$ almost everywhere.

\textit{Step 2: Equality in \eqref{eq:nonlinear} forces $v=\sqrt{u_1^2+u_2^2}$ to be an optimizer of the single GNS inequality.} Indeed, equality in \eqref{eq:nonlinear} requires that $v$ saturates \eqref{lem:single_gns} and that equality holds in the convexity inequality \eqref{convex}, namely
\begin{equation}
\norm{(-\Delta)^{\frac s2} v}{2}^2 = \norm{(-\Delta)^{\frac s2} u_1}{2}^2 + \norm{(-\Delta)^{\frac s2} u_2}{2}^2.
\end{equation}
Consequently, $v(x)=\beta Q(\alpha x+\eta)$ for some $\beta>0$, $\alpha>0$ and $\eta\in\R^N$.

\textit{Step 3: Form of $u_1$ and $u_2$.} From Step 1, $\abs{u_2}=c\abs{u_1}$, so $v = \sqrt{1+c^2} \abs{u_1}$. Thus
\begin{equation}
\abs{u_1} = \frac{\beta}{\sqrt{1+c^2}} Q(\alpha x+\eta), \quad \abs{u_2} = \frac{c\beta}{\sqrt{1+c^2}} Q(\alpha x+\eta).
\end{equation}
Setting $\kappa=\beta$ and defining $\theta$ by $\sin\theta = 1/\sqrt{1+c^2}$, $\cos\theta = c/\sqrt{1+c^2}$ gives the desired form. The normalization condition follows from a direct computation.
\end{proof}
\subsection*{Conflict of Interest}
The authors declare no conflict of interest.

\subsection*{Data Availability Statement}
Not applicable. No data were generated or analyzed in this study.


\begin{thebibliography}{99}
	
	\bibitem{BartschJeanjeanSoave2016}
	\href{https://www.sciencedirect.com/science/article/pii/S0021782416300046}{T.~Bartsch, L.~Jeanjean and N.~Soave,
		\newblock Normalized solutions for a system of coupled cubic Schr\"odinger equations on $\mathbb R^3$,
		\newblock J. Math. Pures Appl. \textbf{106} (2016), 583--614.}
	
	\bibitem{Dalfovo1999}
	\href{https://doi.org/10.1103/RevModPhys.71.463}{F.~Dalfovo, S.~Giorgini, L.~P. Pitaevskii and S.~Stringari,
		\newblock Theory of Bose--Einstein condensation in trapped gases,
		\newblock Rev. Mod. Phys. \textbf{71} (1999), 463.}
	
	\bibitem{DiNezzaPalatucciValdinoci2012}
	\href{https://www.sciencedirect.com/science/article/pii/S0007449711001254}{E.~Di Nezza, G.~Palatucci and E.~Valdinoci,
		\newblock Hitchhiker's guide to the fractional Sobolev spaces,
		\newblock Bull. Sci. Math. \textbf{136} (2012), 521--573.}
	
	\bibitem{DuTianWangZhang2019}
	\href{https://doi.org/10.1017/prm.2018.41}{M.~Du, L.~Tian, J.~Wang and F.~Zhang,
		\newblock Existence of normalized solutions for nonlinear fractional Schr\"odinger equations with trapping potentials,
		\newblock Proc. Roy. Soc. Edinburgh Sect. A \textbf{149} (2019), 617--653.}
	
	\bibitem{FrankLenzmann2013}
	\href{https://doi.org/10.1007/s11511-013-0095-9}{R.~L. Frank and E.~Lenzmann,
		\newblock Uniqueness of non-linear ground states for fractional Laplacians in $\mathbb R$,
		\newblock Acta Math. \textbf{210} (2013), 261--318.}
	
	\bibitem{FrankLenzmannSilvestre2016}
	\href{https://doi.org/10.1002/cpa.21591}{R.~L. Frank, E.~Lenzmann and L.~Silvestre,
		\newblock Uniqueness of radial solutions for the fractional Laplacian,
		\newblock Comm. Pure Appl. Math. \textbf{69} (2016), 1671--1726.}
	
	\bibitem{GuoLiWeiZeng2019}
	\href{https://doi.org/10.1016/j.jfa.2018.09.015}{Y.~J. Guo, S.~Li, J.~C. Wei and X.~Y. Zeng,
		\newblock Ground states of two-component attractive Bose--Einstein condensates I: Existence and uniqueness,
		\newblock J. Funct. Anal. \textbf{276} (2019), 183--230.}
	
	\bibitem{GuoSeiringer2014}
	\href{https://doi.org/10.1007/s11005-013-0667-9}{Y.~J. Guo and R.~Seiringer,
		\newblock On the mass concentration for Bose--Einstein condensates with attractive interactions,
		\newblock Lett. Math. Phys. \textbf{104} (2014), 141--156.}
	
	\bibitem{GuoWangZengZhou2018}
	\href{https://iopscience.iop.org/article/10.1088/1361-6544/aa99a8}{Y.~J. Guo, Z.~Q. Wang, X.~Y. Zeng and H.~S. Zhou,
		\newblock Properties of ground states of attractive Gross--Pitaevskii equations with multi-well potentials,
		\newblock Nonlinearity \textbf{31} (2018), 957--979.}
	
	\bibitem{GuoZengZhou2016}
	\href{https://www.sciencedirect.com/science/article/pii/S0294144915000165}{Y.~J. Guo, X.~Y. Zeng and H.~S. Zhou,
		\newblock Energy estimates and symmetry breaking in attractive Bose--Einstein condensates with ring-shaped potentials,
		\newblock Ann. Inst. H. Poincar\'e Anal. Non Lin\'eaire \textbf{33} (2016), 809--828.}
	
	\bibitem{GuoZengZhou2017}
	\href{https://www.aimsciences.org/article/doi/10.3934/dcds.2017159}{Y.~J. Guo, X.~Y. Zeng and H.~S. Zhou,
		\newblock Blow-up solutions for two coupled Gross--Pitaevskii equations with attractive interactions,
		\newblock Discrete Contin. Dyn. Syst. \textbf{37} (2017), 3749--3786.}
	
	\bibitem{GuoZengZhou2018}
	\href{https://doi.org/10.1016/j.jde.2017.09.039}{Y.~J. Guo, X.~Y. Zeng and H.~S. Zhou,
		\newblock Blow-up behavior of ground states for a nonlinear Schr\"odinger system with attractive and repulsive interactions,
		\newblock J. Differential Equations \textbf{264} (2018), 1411--1441.}
	
	\bibitem{Ho1996}
	\href{https://doi.org/10.1103/PhysRevLett.77.3276}{T.-L. Ho and V.~B. Shenoy,
		\newblock Binary mixtures of Bose condensates of alkali atoms,
		\newblock Phys. Rev. Lett. \textbf{77} (1996), 3276--3279.}
	
	\bibitem{Kim2025}
	\href{https://doi.org/10.58997/ejde.2025.83}{Y.-C.~Kim,
		\newblock H\"older regularity of weak solutions to nonlocal \(p\)-Laplacian type Schr\"odinger equations with \(A_{1}^{p}\)-Muckenhoupt potentials,
		\newblock Electron. J. Differ. Equ. \textbf{2025} (2025), No.~83, 1--26.}
	
	\bibitem{Laskin2000}
	\href{https://doi.org/10.1016/S0375-9601(00)00201-2}{N.~Laskin,
		\newblock Fractional quantum mechanics and L\'evy path integrals,
		\newblock Phys. Lett. A \textbf{268} (2000), 298--305.}
	
	\bibitem{Laskin2002}
	\href{https://doi.org/10.1103/PhysRevE.66.056108}{N.~Laskin,
		\newblock Fractional Schr\"odinger equation,
		\newblock Phys. Rev. E \textbf{66} (2002), 056108.}
	
	\bibitem{LiuZhangZuoCPAA}
	\href{https://www.aimsciences.org/article/doi/10.3934/cpaa.2025061}{C.~Liu, Z.~Zhang and J.~Zuo,
		\newblock Existence and non-existence of normalized solutions for a nonlinear fractional Schr\"odinger system,
		\newblock Commun. Pure Appl. Anal. \textbf{24} (2025), no.~10, 1830--1860.}
	
	\bibitem{LiuZhangZuoUniqueness}
	\href{https://doi.org/10.1007/s00245-026-10507-y}{C.~Liu, Z.~Zhang and J.~Zuo,
		\newblock The uniqueness and concentration behavior of solutions for a nonlinear fractional Schr\"odinger system,
		\newblock Appl. Math. Optim. \textbf{94} (2026), 56.}
	
	\bibitem{LuoZhang2020}
	\href{https://doi.org/10.1007/s00526-020-01814-5}{H.~Luo and Z.~Zhang,
		\newblock Normalized solutions to the fractional Schr\"odinger equations with combined nonlinearities,
		\newblock Calc. Var. Partial Differential Equations \textbf{59} (2020), 35 pp.}
	
	\bibitem{PengXia2021}
	\href{https://www.aimsciences.org//article/doi/10.3934/cpaa.2021128}{S.~Peng and A.~Xia,
		\newblock Normalized solutions of supercritical nonlinear fractional Schr\"odinger equation with potential,
		\newblock Commun. Pure Appl. Anal. \textbf{20} (2021), 3723--3744.}
	
	\bibitem{RosOtonSerra2014}
	\href{https://doi.org/10.1007/s00205-014-0740-2}{X.~Ros-Oton and J.~Serra,
		\newblock The Pohozaev identity for the fractional Laplacian,
		\newblock Arch. Rational Mech. Anal. \textbf{213} (2014), 587--628.}	
	
	\bibitem{Secchi2016}
	\href{https://www.projecteuclid.org/journals/topological-methods-in-nonlinear-analysis/volume-47/issue-1/On-fractional-Schr%c3%b6dinger-equations-in-mathbbRN-without-the-Ambrosetti-Rabinowitz/10.12775/TMNA.2015.090.full}{S.~Secchi,
		\newblock On fractional Schr\"odinger equations in \(\mathbb{R}^N\) without the Ambrosetti--Rabinowitz condition,
		\newblock Topol. Methods Nonlinear Anal. \textbf{47} (2016), no.~1, 19--41.}	
\end{thebibliography}
\end{document}